\documentclass[10pt]{article}

\usepackage{bm}
\usepackage{amsmath}
\usepackage{amsthm}
\usepackage{amsfonts}
\usepackage{amssymb}
\usepackage{amscd}
\usepackage{mwe}
\usepackage{cite}
\usepackage{multicol}
\usepackage{hyperref}
\usepackage{pstricks}
\usepackage{latexsym}
\usepackage{enumerate}
\usepackage{epsfig}
\usepackage{subfig}
\usepackage{graphicx}
\usepackage{caption}
\usepackage{framed,fancybox}
\usepackage{bbm}
\usepackage{multirow}
\usepackage[T1]{fontenc}
\usepackage{threeparttable}
\usepackage{rotating}
\usepackage{color,soul}
\usepackage[version=4]{mhchem}
\usepackage{siunitx}
\usepackage{footnpag}			      	
\usepackage{longtable,tabularx}
\newcolumntype{C}{ >{\centering\arraybackslash} m{4cm} }
\newcolumntype{D}{ >{\centering\arraybackslash} m{1cm} }
\newcommand{\R}{\mathbb R}

\newcommand{\crb}[1]{\{ #1 \}}

\newcommand{\m}[1]{{\bf{#1}}}

\newcommand{\tpose}{^{\sf T}}
\newcommand{\C}[1]{{\cal {#1}}}

\newcommand{\J}{\mathcal{J}}
\newcommand{\dt}[1]{\dfrac{\text{d} #1}{\text{d}t}}
\newcommand{\M}{\mathcal{M}}
\newcommand{\deriv}[2]{\dfrac{\text{d} #1}{\text{d} #2}}

\allowdisplaybreaks

\makeatletter
\newcommand{\thickhline}{%
    \noalign {\ifnum 0=`}\fi \hrule height 1pt
    \futurelet \reserved@a \@xhline
}
\newcolumntype{"}{@{\hskip\tabcolsep\vrule width 1pt\hskip\tabcolsep}}
\makeatother

\definecolor{shadecolor}{gray}{0.99}
{\endMakeFramed}

\long\def\symbolfootnote[#1]#2{\begingroup\def\thefootnote{\fnsymbol{footnote}}\footnote[#1]{#2}\endgroup}
\def\qed{\hfill{$\vcenter{\hrule height1pt \hbox{\vrule width1pt height5pt
    \kern5pt \vrule width1pt} \hrule height1pt}$} \medskip}
\title{\bf Numerical Optimization Study of a Constrained Hypersonic Reentry Vehicle}

\author{Cale A.~Byczkowski\thanks{Ph.D.~Candidate, Department of Mechanical and Aerospace Engineering, University of Florida, Gainesville, Florida 32611-6250. Email: cale.byczkowski@ufl.edu.} \\
Anil V.~Rao\thanks{Professor, Department of Mechanical and Aerospace Engineering, University of Florida, Gainesville, FL 32611-6250. E-mail:  anilvrao@ufl.edu.  Corresponding Author.} \vspace{12pt} \\ {\em University of Florida} \\ {\em Gainesville, FL 32611}
}

\date{}

\begin{document}
\maketitle
\thispagestyle{empty}

\begin{abstract}
The trajectory optimization of the atmospheric entry of a reusable launch vehicle is studied.  The objective is to maximize the crossrange of the vehicle subject to two control-inequality path constraints, two state-inequality path constraints, and one mixed state-and-control inequality path constraint.  In order to determine the complex switching structure in the activity of the path constraints, a recently developed method for solving state-path constrained optimal control problems is used. This recently developed method is designed to algorithmically locate the points of activation and deactivation in the path constraints and partition the domain of the independent variable into subdomains based on these activation and deactivation points. Additionally, in a domain where a state-inequality path constraint is found to be active, the method algorithmically determines and enforces the additional necessary conditions that apply on the constrained arc.  A multiple-domain formulation of Legendre-Gauss-Radau direct collocation is then employed to transcribe the optimal control problem into a large sparse nonlinear programming problem. Two studies are performed which analyze a variety of problem formulations of the hypersonic reusable launch vehicle. Key features of the constrained trajectories are presented, and the method used is shown to obtain highly accurate solutions with minimal user intervention.

\end{abstract}


\renewcommand{\baselinestretch}{1}
\normalsize\normalfont 

\section{Nomenclature}

{\renewcommand\arraystretch{1.0}
 \noindent
 \begin{longtable*}{@{}l @{\quad=\quad} l@{}}
 $A$ & vehicle reference area, m$^2$ \\
 $\bf{b}$ & boundary conditions \\
 $\bf{c}$ & inequality path constraints \\
 $C_D$ & coefficient of drag \\
 $C_L$ & coefficient of lift \\
 $d$ & domain index \\
 $\bf{D}$ & Legendre-Gauss-Radau differentiation matrix \\
 $D$ &  total number of domains \\
$\bf{f}$ & right-hand side of dynamics \\
 $h$ & altitude over spherical Earth, m \\
 $H_s$   & scale height, m \\
 $\C{H}$   & Hamiltonian of optimal control problem \\
 $i$ & collocation point row index \\
 $j$  & collocation point column index \\
 $\C{J}$  & objective functional \\
 $k$ & mesh index \\
 $K^{\crb{d}}$ & number of mesh intervals in a domain $d$ \\
 $\ell_{j}^{(k)}$ & Lagrange polynomial $j$ of mesh interval $k$ \\
 $\C{L}$& integrand appearing in Lagrange cost \\
 $\C{M}$ & Mayer cost \\
 $M$ & mesh iteration count \\
 $m$ & vehicle mass \\
 $n$ & sensed acceleration load \\
 $n_c$ & number of inequality path constraints \\
 $n_{sc}$ & number of state-variable inequality path constraints \\
 $n_u$ & number of controls \\
 $n_y$ & number of states \\
 $N^{\crb{d}}$ & total number of collocation points in domain $d$ \\
 $\m{N}$ & tangency constraint vector \\
 $\C{P}_{d}$ & time domain $d$ \\
 $q$ & dynamic pressure, kPa, or order of state-variable inequality path constraint \\
 $q_{\max}$ & maximum dynamic pressure, kPa \\
 $\dot{Q}$ & stagnation point heating rate, W/m$^2$ \\
 $\hat{\dot{Q}}$ & stagnation point heating rate multiplier, W/m$^2$ \\ 
 $\dot{Q}_{\max}$ & maximum allowable stagnation point heating rate, W/m$^2$ \\
 $r$ & radial distance from the center of Earth, m \\
 $r_n$ & nose radius, m \\
 $R_e$ & radius of Earth, m \\
 $s_i$ & $i$-th state-path constraint \\
 $\mathcal{S}_{k}$ & mesh interval $k$ \\
 $t$ & time on $[t_0,t_f]$ \\
 $t_s^{\crb{d}}$ & domain interface variable \\
 $\m{u}$ & control in time horizon $t \in [t_0,t_f]$\\
 $\m{U}$ & matrix of control characterization at LGR collocation points \\
 $v$ & planet-relative speed, m/s \\
 $\m{w}$ & LGR weights \\
 $\m{y}$ & state in time horizon $t \in [t_0,t_f]$ \\
 $\m{Y}$ & matrix of state approximation at LGR points \\
 $\alpha$ & angle of attack, deg or rad \\
 $\beta$ & inverse of density scale height, 1/m \\
 $\gamma$ & planet-relative flight path angle, deg or rad\\
 $\epsilon$ & detection tolerance \\
 $\epsilon_{mesh}$ & mesh error parameter \\
 $\epsilon_{c}$ & constraint error parameter \\
 $\theta$ & planet-relative longitude, deg or rad \\
 $\mu$ & gravitational parameter, m$^3$/s$^2$ \\
 $\rho$ & atmospheric density, kg/m$^3$ \\
 $\rho_0$ & atmospheric density at sea level, kg/m$^3$ \\
 $\tau$ & time on $[-1,+1]$ \\
 $\sigma$ & bank angle, deg or rad \\
 $\phi$ & planet-centered latitude, deg or rad \\
 $\psi$ & azimuth, deg or rad \\
 \end{longtable*}}

\renewcommand{\baselinestretch}{1.5}
\normalsize\normalfont
\section{Introduction}
Trajectory optimization of a reusable launch vehicle entering the atmosphere of the Earth dates back to the Space Shuttle program \cite{Dickmanns1972,ZondervanHuffman1984}. Since that time, the reusable launch vehicle entry (RLVE) problem has served as a challenging real-world application for researchers in the field of optimal control \cite{KugelmannPesch1990,Pesch1994,Betts2020}. Although many solutions to the RLVE problem have been obtained, the RLVE problem remains an excellent benchmark problem for newly developed numerical methods for optimal control \cite{Betts2020,DarbyRao2011b,PattersonRao2015,DennisRao2019,MallTaheri2022}. Furthermore, key solution characteristics of the RLVE problem, (for example, segments where inequality path constraints are active) appear in a variety of other applications including aeroassisted orbital transfer \cite{RaoTang2002,DarbyRao2011c,FuhrRao2018}, hypersonic problems for both Earth and Mars entry \cite{GrantBraun2015,ZhengAi2017}, and boost-glide vehicles \cite{MillerRao2021}. Depending upon the formulation of the RLVE problem, the solution may become increasingly more challenging to compute. First, the RLVE problem may already be challenging to solve even if the problem does not contain any inequality path constraints because the dynamics are coupled and highly nonlinear. Second, this challenge increases when an RLVE problem is subject to inequality constraints on the control (for example, limits on the bank angle and angle of attack).  Third, the challenge of solving an RLVE problem increases even further when the problem is subject to inequality path constraints that are functions of the state (for example, heating rate, dynamic pressure, and sensed acceleration constraints). The aforementioned challenges, including the increasing challenge as additional constraints are enforced, necessitates the need to develop advanced numerical methods for optimal control in order to determine accurate solutions to problems of this kind.

Numerical methods for solving optimal control problems fall into the following two broad categories:  indirect methods and direct methods. In an indirect method, the calculus of variations is applied which leads to a Hamiltonian boundary value problem (HBVP) that is solved numerically \cite{BrysonHo1975}. In a direct method, the optimal control problem is approximated as a finite-dimensional nonlinear programming problem (NLP) using either a control or a state and control parameterization \cite{Betts1998}. Direct methods have several advantages over indirect methods such as having a larger domain of convergence and allowing for flexibility in the problem formulation when including path constraints.

Concurrent with the advancement of direct methods for solving challenging optimal control problems, a great amount of research has been conducted over the past decade directed towards developing new techniques for improving the tractability of indirect methods \cite{BertrandEpenoy2002,GraichenChaplais2010}. Indirect methods have the advantage that any solution to the resulting HBVP is guaranteed to be a local extremal. As discussed, however, in the existing literature, indirect methods have two major drawbacks. First, indirect methods have a small radius of convergence due to the need to supply an accurate initial guess of the boundary values of the costate. Second, indirect methods are challenging to implement in the presence of path constraints due to the required knowledge of the switching structure in the path constraint activity/inactivity. Not only is this switching structure generally not known a priori, but imposing any switching structure as part of the indirect method makes an original two-point boundary-value problem (TPBVP) a multiple-point boundary-value problem (MPBVP). Moreover, this MPBVP changes depending upon the number and location of active inequality constraints (which, as stated, is not known a priori) thus making indirect methods difficult to implement. 

To circumvent both major drawbacks of indirect methods, recent research has been aimed at developing advanced numerical continuation/homotopy methods \cite{Epenoy2011,GrantBraun2015,TaheriAtkins2016,ZhengAi2017,AntonyGrant2018,MansellGrant2018,MallTaheri2020,HeidrichGrant2021}. A homotopy approach is based on the idea of solving a series of simpler problems and successively approaching the desired complex problem of interest by varying a set of continuation/homotopy parameters. While the use of continuation methods have been shown to make indirect methods more tractable when solving challenging trajectory optimization problems, there still requires a large amount of user input and attention to detail. For example, Ref. \cite{MallTaheri2022} adapts the recently developed uniform trigonometrization method (UTM) \cite{MallTaheri2020} to solve a complex constrained shuttle-type entry vehicle trajectory optimization problem. In order to generate feasible solutions and accurately capture the optimal constraint structures, the UTM performs a large set of continuation steps requiring the solution of an even larger number of TPBVPs.

An alternative to indirect methods is the class of direct collocation methods \cite{Betts2020}.  In a direct collocation method, the state and control are parameterized using trial or basis functions and the dynamics are approximated (collocated) at a finite set of points (nodes).  While many collocation methods have been developed, over the past few decades the class of {\em Gaussian quadrature orthogonal direct collocation methods} \cite{BensonRao2006,KameswaranBiegler2008a,GargHuntington2010,GargRao2011a,GargRao2011b,HagerWang2018} has been developed. The most well developed Gaussian quadrature collocation methods use either Legendre-Gauss (LG) points~\cite{BensonRao2006}, Legendre-Gauss-Radau (LGR) points~\cite{KameswaranBiegler2008a,GargHuntington2010,GargRao2011a,GargRao2011b}, or Legendre-Gauss-Lobatto (LGL) points~\cite{Elnagar1995}. Although these methods converge exponentially when the solution of the optimal control problem satisfies particular smoothness and coercivity conditions~\cite{HouRao2012}, they often struggle to compute an accurate solution when the underlying solution is nonsmooth. To account for this, a class of so-called $hp$-adaptive mesh refinement techniques have been developed to exploit the solution structure and improve solution accuracy\cite{DarbyRao2011a,DarbyRao2011b,PattersonRao2015,LiuRao2015,LiuRao2018}. While $hp$-adaptive methods are highly accurate relative to traditional fixed-order methods, these methods may not produce sufficiently accurate solutions for some classes of problems. Specifically, when the optimal control problem includes control and/or state and control constraints, the resulting NLP may become ill-conditioned, leading to control chattering behavior \cite{Betts2020,PagerRao2022,ByczkowskiRao2024a}. 

The objective of this paper is to apply a newly developed adaptation of the Legendre-Gauss-Radau direct collocation method for solving state-path constrained optimal control problems.  This newly developed method, called the SPOC method \cite{ByczkowskiRao2024a}, handles the state-path constraints algorithmically by enforcing the so-called {\em tangency conditions} \cite{BrysonHo1975} to perform index reduction on the resulting high-index system of differential algebraic equations \cite{FeeheryBarton1998}. Furthermore, the SPOC method employs a structure detection algorithm to decompose the original time domain into subdomains consisting of active and inactive constraint subarcs \cite{ByczkowskiRao2023a}.  The combination of algorithmic index reduction together with domain decomposition is then used to formulate a multiple-domain optimal control problem where the domains reflect segments where path constraints are either active or inactive.  This multiple-domain formulation also introduces the start and termination times of the path constraint activity as decision variables.  The multiple-domain formulation of the optimal control problem is then approximated using Legendre-Gauss-Radau collocation and is subsequently transcribed into a large sparse NLP.  

The contributions of this work are as follows. First, this research is the first instance of applying the SPOC method to a problem with multiple active state-path constraints throughout the trajectory. Furthermore, this study is the first instance of the SPOC method solving a problem with control-inequality path constraints and mixed state-control inequality path constraints present in the problem formulation. Second, two studies are performed that solve a variety of formulations of the RLVE problem. Insights are provided into the change of the active state-path constraint structure when (i) control-inequality path constraints are enforced and (ii) the rotation of the Earth is included. The results display the ability of the SPOC method to algorithmically detect and optimize the constraint structures without a priori information of the solution structure. Third, the results obtained using the SPOC method approximate the active state-path constraints to a higher resolution than current state-of-the-art direct methods. Moreover, numerical comparisons are made against the results provided in Ref. \cite{MallTaheri2022} to further validate the accuracy of the SPOC solution. The fourth contribution of this work is to highlight a feature of the SPOC method not discussed in previous work, which is the ability to capture small control discontinuities due to state-path constraints becoming active/inactive. The multiple-domain LGR formulation used within the SPOC method allows for the computation of separate control arcs leading to an improved approximation of the active state-path constraints. Finally, it is noted that a preliminary version of the research presented in this paper is described in Ref.~\cite{ByczkowskiRao2024b}.

The remainder of the paper is organized as follows. Section~\ref{section:overview} provides a brief overview of the method for solving state-constrained optimal control problems \cite{ByczkowskiRao2024a}.  Section~\ref{section:probform} formulates the aforementioned hypersonic RLVE trajectory optimization problem. Section~\ref{section:results} demonstrates the SPOC method for a variety of RLVE problem formulations and compares the solutions obtained using the SPOC method against solutions obtained using both an $hp$-adaptive method~\cite{LiuRao2018} and an advanced indirect method~\cite{MallTaheri2022}. Finally, Section~\ref{section:conclusions} provides conclusions on this research.

\section{Method for State-Path Constrained Optimal Control}\label{section:overview}
In this section, a previously developed method for solving state-path constrained optimal control problems \cite{ByczkowskiRao2024a} (SPOC) is reviewed. First, Section~\ref{section:Bolza} presents a general optimal control problem in Bolza form. Next, the multiple-domain Legendre-Gauss-Radau (LGR) direct collocation scheme used to transcribe the continuous time Bolza optimal control problem into a large sparse nonlinear programming problem (NLP) is provided in Section~\ref{section:mdlgr}. The remaining subsections highlight the stages of the SPOC method that reformulate the NLP to include additional constraints and decision variables for solving problems involving state-path constraints. For further details on the SPOC method, the reader is referred to Ref.~\cite{ByczkowskiRao2024a}.

\subsection{Bolza Optimal Control Problem\label{section:Bolza}}
Consider the following general optimal control problem in Bolza form. Determine the state, $\m{y}(t) \in \R^{n_y}$, the control, $\m{u}(t) \in \R^{n_u}$, the initial time, $t_0$, and the terminal time, $t_f$, on the time interval, $t \in [t_0, t_f]$, (where $n_y$ and $n_u$ represent the number of states and controls, respectively) that minimizes the cost functional
\begin{equation}\label{eq:cost1}
\J = \M(\m{y}(t_0),t_0,\m{y}(t_f),t_f) + \int_{t_0}^{t_f} \mathcal{L}(\m{y}(t),\m{u}(t),t) \text{d} t,
\end{equation}
subject to the dynamic constraints
\begin{equation}\label{eq:optprb1_dyn}
\dt{\m{y}} = \m{f}(\m{y}(t),\m{u}(t),t),
\end{equation}
the inequality path constraints 
\begin{equation}\label{eq:optprb1_ineq}
\m{c}_{\text{min}} \leq \m{c}(\m{y}(t),\m{u}(t),t) \leq \m{c}_{\text{max}}, 
\end{equation}
and the boundary conditions
\begin{equation}\label{eq:bc1}
\m{b}_{\text{min}} \leq \m{b}(\m{y}(t_0),t_0,\m{y}(t_f),t_f) \leq \m{b}_{\text{max}}.
\end{equation}
\subsection{Multiple-Domain Legendre-Gauss-Radau Collocation\label{section:mdlgr}}
While any valid collocation scheme can be used, the SPOC method uses Legendre-Gauss-Radau (LGR) collocation due to its ability to converge to a solution at an exponential rate when the underlying solution is smooth \cite{HouRao2012}. However, for the problem studied in this work, the solution may be nonsmooth due to state-path constraints becoming active. Consequently, in order to enforce additional necessary conditions and handle potential nonsmooth behavior, the SPOC method utilizes a multiple-domain LGR collocation scheme described as follows. 

Let the time horizon, $t \in [t_0,t_f]$ of the previously defined Bolza optimal control problem given in Eqs.~\eqref{eq:cost1}--\eqref{eq:bc1} be partitioned into $D$ domains, $\C{P}_d = [t_s^{\crb{d-1}},t_s^{\crb{d}}] \subseteq [t_0,t_f]$, $d = 1,\ldots,D$, such that
\begin{equation}
  \begin{array}{cc}
    \displaystyle \bigcup_{q=1}^{D} \C{P}_q = [t_0,t_f],& \displaystyle \bigcap_{q=d}^{d+1} \C{P}_q = \crb{t_s^{\crb{d}}
    },
  \end{array}
\end{equation}
where $t_s^{\crb{d}}$, for $d = 1,\ldots,D-1$, are the \textit{domain interface variables}, $t_s^{\crb{0}}=t_0$ and $t_s^{\crb{D}}=t_f$, are the initial and final times, respectively, $d$ is the domain index, and the subscript $s$ signifies a switch time between a constrained or unconstrained domain. The domain interface variables $t_s^{\crb{d}}$, for $d = 1,\ldots,D-1$, become additional decision variables in the reformulated NLP and are \textit{not} treated as collocation points. Specifically, the domain interface variables represent either an entry or exit point of an active state-path constraint. \par
Each domain, $\C{P}_d = [t_s^{\crb{d-1}},t_s^{\crb{d}}]$, is then mapped to $\tau \in [-1,+1]$ using the affine transformation in terms of the domain interface variables
\begin{equation}\label{eq:affine_map}
t = \dfrac{t_s^{\crb{d}} - t_s^{\crb{d-1}}}{2} \tau + \dfrac{t_s^{\crb{d}} - t_s^{\crb{d-1}}}{2}.
\end{equation}
Next, the time interval $\tau \in [-1,+1]$ on each domain is divided into $K$ mesh intervals, $\mathcal{S}_k = [T_{k-1},T_{k}] \subseteq [-1,+1]$, for $k = 1,\ldots K$, such that
\begin{equation}
\begin{array}{cc}
  \displaystyle \bigcup_{q=1}^{K} \mathcal{S}_q = [-1,+1], & \displaystyle \bigcap_{q=k}^{k+1} \mathcal{S}_q = \crb{T_k}, \quad (k = 1,\ldots,K-1),
\end{array}
\end{equation}
and $-1 = T_0 < T_1 < \cdots < T_{K-1} < T_K = +1$. Once the multiple-domain formulation is complete, LGR collocation is used to discretize the continuous time problem on each mesh interval. For each mesh interval, the LGR points are defined on $[T_{k-1},T_{k}] \subseteq [-1,+1]$, $k = 1,\ldots,K$. The state of the continuous time optimal control problem is then approximated within each mesh interval $\C{S}_k$ by
\begin{equation}\label{eq:state_approx}
\m{y}^{(k)}(\tau) \approx \m{Y}^{(k)}(\tau) = \sum_{j=1}^{N_k + 1}\m{Y}_j^{(k)}\ell_j^{(k)}(\tau), \quad \ell_j^{(k)}(\tau) = \prod_{\substack{i = 1 \\ i \neq j}}^{N_k+1} \dfrac{\tau - \tau_i^{(k)}}{\tau_j^{(k)}-\tau_i^{(k)}}, \quad (k = 1,\ldots,K),
\end{equation}
where $\ell_j^{(k)}(\tau)$, $j = 1,\ldots,N_k + 1$, is a basis of Lagrange polynomials on $S_k$, $(\tau_1^{(k)},\ldots,\tau_{N_k}^{(k)})$ are the set of $N_k$ LGR collocation points in the interval $[T_{k-1},T_k)$, and $\tau^{(k)}_{N_k+1} = T_k$ is a non-collocated support point. Note, the first collocation point of the ($k$+1)-th interval, $\tau_1^{(k+1)} = T_k$, is equal to the non-collocated point in the $k$-th interval, $\tau_{N_k + 1}^{(k)} = T_{k}$. Differentiating Eq.~\eqref{eq:state_approx} with respect to $\tau$ leads to
\begin{equation}
\deriv{\m{Y}^{(k)}(\tau)}{\tau} = \sum_{j=1}^{N_k+1} \m{Y}_j^{(k)} \deriv{\ell_j^{(k)}(\tau)}{\tau}.
\end{equation}
The dynamic constraints given in Eq.~\eqref{eq:optprb1_dyn} are then approximated at the $N_k$-LGR points in mesh interval $\C{S}_k$, $k = 1,\ldots,K$, on domain $ \C{P}_d, d = 1,\ldots,D$, by
\begin{equation}\label{approxdynconstr}
\sum_{j=1}^{N_k+1}D_{ij}^{(k)}\m{Y}_j^{(k)} - \dfrac{t_s^{\crb{d}} - t_s^{\crb{d-1}}}{2} \m{f}\left( \m{Y}_i^{(k)}, \m{U}_i^{(k)}, t(\tau_i^{(k)},t_s^{\crb{d-1}},t_s^{\crb{d}}) \right) = 0, \hspace{0.5em} \quad (i =1,\ldots,N_k)
\end{equation}
where
\begin{equation}
D_{ij}^{(k)} = \dfrac{\text{d} \ell_j^{(k)}(\tau)}{\text{d} \tau}, \hspace{0.5em} \quad (i =1,\ldots,N_k, \quad j = 1,\ldots,N_k+1),
\end{equation}
are the elements of the $N_k \times (N_k + 1)$ Legendre-Gauss-Radau \textit{differentiation matrix} in the mesh interval $\C{S}_k$, and $\m{U}_i^{(k)}$ is the control parameterization at the $i$-th collocation point in mesh interval $\C{S}_k$. Continuity in the state across adjacent mesh intervals $\{ \mathcal{S}_{k-1},\mathcal{S}_{k} \}$, $k = 2,\ldots,K$, and adjacent time domains $\{\C{P}_{d-1},\C{P}_{d}\}$, $d = 2,\ldots,D$, is achieved by using the same variable to represent $\m{Y}_{N_{k-1}+1}^{(k-1)} = \m{Y}^{(k)}_1$, $k = 2,\ldots,K$, and $\m{Y}_{N_{d-1}+1}^{\crb{d-1}} = \m{Y}^{\crb{d}}_1$, $d = 2,\ldots,D$, respectively. \par 
The LGR discretization of the continuous time multiple-domain Bolza optimal control problem results in the following nonlinear programming (NLP). Minimize the objective function
\begin{equation}\label{eq:NLP_cost}
\C{J} = \C{M} \left( \m{Y}_1^{\crb{1}}, t_0, \m{Y}_{N^{\crb{D}}+1}^{\crb{D}}, t_f \right) + \sum_{d=1}^{D}\dfrac{t_s^{\crb{d}}-t_s^{\crb{d-1}}}{2} \left[\m{w}^{\crb{d}}\right]\tpose \m{L}^{\crb{d}},
\end{equation} 
subject to the dynamic constraints
\begin{equation}\label{eq:NLP_dynconstr}
\m{D}^{\crb{d}}\m{Y}^{\crb{d}} - \dfrac{t_s^{\crb{d}}-t_s^{\crb{d-1}}}{2} \m{F}^{\crb{d}} = \m{0}\hspace{0.5em}, \quad (d = 1,\ldots,D),
\end{equation}
the inequality path constraints
\begin{equation}\label{eq:NLP_ineq}
\m{c}_{\text{min}} \leq \m{C}^{\crb{d}} \leq \m{c}_{\text{max}}, \hspace{0.5em} \quad (d = 1,\ldots,D),
\end{equation}
the boundary conditions
\begin{equation}
\m{b}_{\text{min}} \leq \m{b}\left( \m{Y}_1^{\crb{1}}, t_0, \m{Y}_{N^{\crb{D}}+1}^{\crb{D}}, t_f \right) \leq  \m{b}_{\text{max}},
\end{equation}
and the continuity conditions
\begin{equation}\label{eq:contconds}
\m{Y}_{N^{\crb{d-1}}+1}^{\crb{d-1}} = \m{Y}^{\crb{d}}_1, \hspace{0.5em}\quad (d = 2,\ldots,D),
\end{equation}
where $\m{D}^{\crb{d}} \in \R^{N^{\crb{d}} \times (N^{\crb{d}}+1)}$ is the LGR differentiation matrix in the time domain $\C{P}_d$, $d = 1,\ldots,D$, $\m{w}^{\crb{d}} \in \R^{N^{\crb{d}} \times 1}$ are the LGR weights associated with the approximation of the integral of the cost function in time domain $\C{P}_{d}$ via Gaussian quadrature, and $N^{\crb{d}}$ is the total number of collocation points in time domain $\C{P}_d$. As previously mentioned, the continuity conditions in Eq.~\ref{eq:contconds} are implicitly enforced by using the same decision variable in the NLP for $\m{Y}_{N_{d-1}+1}^{\crb{d-1}} = \m{Y}^{\crb{d}}_1$. Lastly, the elements $\m{L}^{\crb{d}} \in \R^{N^{\crb{d}}\times 1}, \m{F}^{\crb{d}} \in \R^{N^{\crb{d}}\times n_y}, \m{C}^{\crb{d}} \in \R^{N^{\crb{d}}\times n_c}$, in Eqs.~\eqref{eq:NLP_cost}-\eqref{eq:NLP_ineq}, are the discrete approximations of the elements $\mathcal{L}$, $\m{f}$, and $\m{c}$ in Eqs.~\eqref{eq:cost1}-\eqref{eq:optprb1_ineq}, respectively, at each collocation point within time domain $\C{P}_{d}$, $d = 1,\ldots,D$. 

\subsection{Structure Detection and Decomposition}\label{section:detection_decomposition}
In order to identify the domain interface variables and partition the original domain into multiple subdomains, the SPOC method employs a structure detection method (SDM) developed in Ref.~\cite{ByczkowskiRao2023a}. The first step in the SDM obtains an approximate solution on a single domain $\C{P}_1 = [t_0,t_f]$ using standard Legendre-Gauss-Radau (LGR) collocation on a fixed mesh consisting of $\mathcal{S}_k,$ $k = 1,\ldots,K$, mesh intervals with $N_k$ collocation points within each mesh interval. Next, let $\m{Y}(\tau_j^{(k)})$, $j \in \crb{1,\ldots,N_k+1}$, $k \in \crb{1,\ldots,K}$, be the approximate state solution obtained on the initial domain at the $j$-th collocation point in the $k$-th mesh interval. Furthermore, suppose that one of the inequality path constraints from Eq.~\eqref{eq:optprb1_ineq} is a state-path constraint of the form $c_i(\m{y}(t),\m{u}(t),t) = s(\m{y}(t),t)$ for some $i \in \crb{1,\ldots,n_c}$. Then the relative difference, $\delta s_i$, is computed at every collocation point plus the final non-collocated point of the final mesh interval for state-path constraint $s_i$, for all $i = 1,\ldots,n_{sc}$, in the problem, given by
\begin{equation}\label{eq:svic_rel_diff}
\delta s_i (\tau_j^{(k)}) =  \dfrac{\left| s_i(\m{Y}(\tau_j^{(k)}),\tau_j^{(k)}) - s_{i,\max / \min} \right|}{1 + |s_{i,\max / \min}|}, \hspace{0.5em} (j = 1, \ldots, N_k + 1, \quad k = 1, \ldots, K),
\end{equation}
where $n_{sc} \leq n_c$ is the number of state-path constraints present in the problem, and $s_{i,\max / \min}$ denotes either the maximum or minimum bound of the state-path constraint $s_i$, $i = 1,\ldots,n_{sc}$. Once the relative differences are obtained, each is compared against a user-specified detection tolerance $\epsilon_i$, $i = 1,\ldots,n_{sc}$, to determine points along the solution where the state-path constraints become active/deactive (A/D). To account for inaccuracies in the solution used to perform detection, as well as to prevent overlapping or collapsing of domains while optimizing the location of the A/D times, bounds on the detected A/D times are computed as follows. Suppose the point $\tau_s = \tau_j^{(k)}$ for some $j \in \{1,\ldots,N_k\}$ and $k \in \{1,\ldots,K\}$ is detected to be either an activation or deactivation point, the bounds are then obtained using
\begin{equation}\label{eq:ADbounds}
  \left.
    \begin{array}{lcl}
      \tau_s^-  &=& \tau_j^{(k)} + \nu \left( \tau_{j-1}^{(k)} - \tau_j^{(k)} \right), \\ \\
      \tau_s^+ &=&  \tau_j^{(k)} + \nu \left( \tau_{j+1}^{(k)} - \tau_j^{(k)} \right),
    \end{array} \right\} \left. \begin{array}{ll}  \\ s \in \{1,\ldots,S\}, \\ \,\end{array} \right.
\end{equation}
where $S$ is the total number of detected A/D points, and $\nu > 0$ is a user-specified value to define the size of the allowable search space while optimizing the exact locations of the domain interface variables. Based on the results of the SDM, the SPOC method automatically partitions the previous solution into multiple subdomains and, within subdomains where the state-path constraints are deemed to be active, algorithmically enforces the additional necessary conditions discussed in the next section. 

\subsection{Domain Constraints and Refinement}\label{section:constraints_refinement}
Within domains that have been identified to have an active state-path constraint, the SPOC method algorithmically includes the additional necessary conditions derived in Bryson, et. al.~\cite{BrysonDreyfus1963}, which are given as follows. For simplicity, the discussion to follow assumes a single state-path constraint is present in the problem formulation, but the application of the SPOC method can be extended to problems with multiple state-path constraints, and is a primary focus of the application studied in this paper.

Suppose that on the optimal solution the upper limit of the state-path constraint, $s_i(\m{y}(t),t)$, for some $i \in \crb{1,\ldots,n_c}$, is active over a non-zero interval of time\footnote[1]{The same analysis can be performed if the lower limit of the state-path constraint is active.}, that is
\begin{equation}\label{eq:activeSVICindirectadjoining}
s_i(\m{y}(t),t) - s_{i,\max} = 0, \quad t \in [t_1, t_2] \quad (\text{for some } i \in \crb{1,\ldots,n_c})
\end{equation}
where $t_1$ and $t_2$ are the activation and deactivation times of the state-path constraint, respectively. Given that Eq.~\eqref{eq:activeSVICindirectadjoining} must be satisfied on the optimal solution, all the time derivatives of $s_i(\m{y}(t),t)$ must also be equal to zero. Thus, successive time derivatives of Eq.~\eqref{eq:activeSVICindirectadjoining} are taken, and the dynamic constraints $\m{f}(\m{y}(t),\m{u}(t),t)$ are substituted for $\text{d}\m{y}(t)/\text{d}t$ until an expression explicit in a control variable is obtained. If $q$ successive time derivatives are required, then the state-path constraint given by Eq.~\eqref{eq:activeSVICindirectadjoining} is said to be of $q$th order. Given the control necessary to prevent $s_i(\m{y}(t),t)$ from being violated is obtained from its $q$th time derivative, no finite control will keep the system on the constraint boundary if the path entering the constraint boundary does not meet the following ``tangency'' constraints~\cite{BrysonHo1975}
\begin{equation}\label{eq:tangencyconstraints}
\m{N}(\m{y}(t_1),t_1) \equiv \begin{bmatrix}
s(\m{y},t_1) - s_{i,\max} \\[5pt] 
\deriv{}{t}s(\m{y},t_1) \\[5pt] 
\vdots \\[5pt] 
\deriv{^{(q-1)}}{t^{(q-1)}}s(\m{y},t_1)
\end{bmatrix}
= \m{0}.
\end{equation}
Due to the activation of a state-path constraint, a high-index system of differential algebraic equations (DAEs) arises. The necessary tangency conditions are analogous to defining a set of  \textit{consistent initial conditions} when performing index-reduction of the resulting high-index system of DAEs. Consequently, the state-path constraint given by Eq.~\eqref{eq:activeSVICindirectadjoining} is replaced by the $q$th time derivative which is now a mixed state-control equality constraint
\begin{equation}\label{eq:mixedconstr}
\deriv{^{(q)}}{t^{(q)}}s(\m{y},t) = c_i(\m{y}(t),\m{u}(t)) = 0, \quad t \in [t_1,t_2].
\end{equation}
Specifically, the $i$-th active inequality state-path constraint in Eq.~\eqref{eq:optprb1_ineq} is converted to an equality constraint by replacing the $i$-th element of Eq.~\eqref{eq:optprb1_ineq} with Eq.~\eqref{eq:mixedconstr}.
\par
After solving the reformulated problem, the SPOC method checks the current solution against two user-defined error tolerances. The first error tolerance, $\epsilon_c$, is a measure of maximum allowable constraint violations while the second error tolerance, $\epsilon_{mesh}$, is a measure of the error in the state approximation. If neither error tolerances are satisfied on the current mesh iteration $M$, the SPOC method employs a standard mesh refinement strategy and re-performs structure detection, thus providing the ability to perform structure detection on a solution that has been obtained on a refined mesh. As will be discussed in Section~\ref{section:results}, the ability to re-perform structure detection for the problem studied in this work will be necessary due to the complex constraint activation/deactivation structure.

\section{Problem Formulation\label{section:probform}}
In this section, the hypersonic reusable launch vehicle entry (RLVE) optimal control problem is presented. The formulation follows closely to that presented in Ref.~\cite{MallTaheri2022}, where the objective is to maximize the achieved crossrange while being subject to two control constraints, two state-path constraints, and one mixed state-control constraint. The remainder of this section presents the equations of motion, path constraints enforced on the vehicle, and the resulting optimal control problem. 
\subsection{Equations of Motion}
The equations of motion for a point mass in motion over a spherical nonrotating planet are given in spherical coordinates as \cite{RaoTang2002,Betts2020}
\begin{subequations}\label{eq:EOM}
\begin{align}
\dot{r} &= v \sin \gamma, \label{eq:rdot}\\[5pt]
\dot{\theta} &= \dfrac{v \cos \gamma \sin \psi}{r \cos \phi}, \label{eq:thetadot}\\[5pt]
\dot{\phi} &= \dfrac{v \cos \gamma \cos \psi}{r}, \label{eq:phidot}\\[5pt]
\dot{v} &= -D -g\sin \gamma, \label{eq:vdot}\\[5pt]
\dot{\gamma} &= \dfrac{L \cos \sigma}{v} + \cos \gamma \left( \dfrac{v}{r} - \dfrac{g}{v}  \right), \label{eq:gammadot}\\[5pt]
\dot{\psi} &=  \dfrac{L \sin \sigma}{v \cos \gamma} + \dfrac{v}{r} \cos \gamma \sin \psi \tan \phi,\label{eq:psidot}
\end{align}
\end{subequations}
where $g = \mu/r^2$, $r$ is the radial distance, $\theta$ is the planet-relative longitude, $\phi$ is the planet-centered latitude, $v$ is the planet-relative speed, $\gamma$ is the planet-relative flight path angle, $\psi$ is the azimuth, $\alpha$ is the angle of attack, and $\sigma$ is the bank angle.  Furthermore, 
\begin{equation}
  \begin{array}{lcl}
    D &=& q A C_D / m, \\
    L &=& q A C_L / m, 
  \end{array}
\end{equation}
are the drag and lift specific forces, respectively, $A$ is the vehicle reference area, $q = \rho v^2 / 2$ is the dynamic pressure,
\begin{equation}
  \rho = \rho_0 \exp \left(-h/H_s \right)
\end{equation}
is the atmospheric density, $h = r - R_e$ is the altitude, $C_D$ is the drag coefficient, and $C_L$ is the lift coefficient. Next, the coefficients of lift and drag are modeled, respectively, as \cite{Betts2020},
\begin{equation}\label{eq:density}
\begin{array}{cll}
C_L &=& C_{L0} + C_{L1}\alpha, \\[5pt]
C_D &=& C_{D0} + C_{D1}\alpha + C_{D2}\alpha^2,
\end{array}
\end{equation}
The aerodynamic coefficients and physical constants used in this model are provided in Table~\ref{table:parameters}.

\begin{table}[htb]
\centering
\caption{Physical and aerodynamic constants}
\begin{tabular}{ c c l } 
 \hline \hline
 \rule{0pt}{2.5ex} Quantity & Value & Unit \\ \hline 
 \rule{0pt}{2ex} $R_e$ & 6371.2039 & km \\[3pt] 
 \rule{0pt}{2ex}  $H_s$ & 7254.24 & m  \\[3pt]
 \rule{0pt}{2ex}  $\rho_0$ & 1.2256 & kg/m$^3$  \\[3pt]
 \rule{0pt}{2ex}  $\mu$ & $3.986031954 \times 10^{5}$ & km$^3$/s$^2$  \\[3pt]
 \rule{0pt}{2ex}  $\omega_e$ & $7.292115856 \times 10^{-5}$ & rad/s  \\[3pt]
 \rule{0pt}{2ex}  $g_0$ & 9.8066498 & m/s$^2$  \\[3pt]
 \rule{0pt}{2ex}  $m$ & 92079.2525 & kg  \\[3pt]
 \rule{0pt}{2ex}  $A$ & 249.9092 & m$^2$  \\
 \rule{0pt}{2ex}  $k$ & 1.7415$\times$10$^{-4}$ & kg$^{1/2}$/m$^2$  \\[3pt]
 \rule{0pt}{2ex}  $r_n$ & 1  &  m   \\[3pt]
 \rule{0pt}{2ex}  $C_{L0}$ & $-0.2070$ & --  \\[3pt]
 \rule{0pt}{2ex}  $C_{L1}$ & 1.6756 & 1/rad  \\[3pt]
 \rule{0pt}{2ex}  $C_{D0}$ & $0.0785$ & --  \\[3pt]
 \rule{0pt}{2ex}  $C_{D1}$ & $-0.3529$ & 1/rad \\[3pt]
 \rule{0pt}{2ex}  $C_{D2}$ & $2.0400$ & 1/rad$^2$ \\[3pt]
  \hline \hline
\end{tabular}
\label{table:parameters}
\end{table}

\subsection{Path Constraints}
The following constraints are placed on the bank angle, $\sigma$, and the angle of attack $\alpha$:
\begin{equation}\label{eq:controlconstraints}
  \begin{array}{lclcl}
    \sigma &\geq& \sigma_{\min}, \\[5pt]
    \alpha &\leq&  \alpha_{\max},
  \end{array}
\end{equation}
where, for the case when both control constraints in Eq.~\eqref{eq:controlconstraints} are enforced, $\sigma_{\text{min}} = -75$ deg and $\alpha_{\max} = 19$~deg. Next, the two state-path constraints enforced are
\begin{equation}\label{eq:heatrate}
\dot{Q} = k \sqrt{\dfrac{\rho}{r_n}} v^3 \leq \dot{Q}_{\max}
\end{equation}
and 
\begin{equation}\label{eq:dynamicpressure}
q = \frac{1}{2} \rho v^2 \leq q_{\max}
\end{equation}
where $\dot{Q}$ is the stagnation point heating rate~\cite{SuttonGraves1971}, $q$ is the dynamic pressure, $r_n$ is the nose radius of the entry vehicle, and $k$ is a heating rate constant. The upper limits on the stagnation point heating rate and dynamic pressure are chosen as $\dot{Q}_{\max} = 0.85$ MW/m$^2$, and $q_{\max} = 12.53$ kPa, respectively. Lastly, the mixed state-control constraint 
\begin{equation}\label{eq:sensedaccel}
n = \dfrac{1}{mg_0} \sqrt{L^2 + D^2}\leq n_{\max}
\end{equation} 
is enforced, where $n$ is the sensed acceleration normalized by sea level gravity $g_0$ and is a function of the angle of attack, $\alpha$. For the results presented, the upper limit on the sensed acceleration was chosen to be $n_{\max} = 1.15$.

\subsection{Optimal Control Problem}

The hypersonic RLVE optimal control problem is as follows. Maximize the final latitude (crossrange)
\begin{equation}\label{eq:objective}
\C{J} = \phi (t_f),
\end{equation} 
subject to the dynamic equation constraints given by Eqs.~\eqref{eq:rdot}-\eqref{eq:psidot}, the control constraints given by Eq.~\eqref{eq:controlconstraints}, the two state-path constraints given by Eqs.~\eqref{eq:heatrate} and.~\eqref{eq:dynamicpressure}, the mixed state-control constraint given by Eq.~\eqref{eq:sensedaccel}, and the following boundary conditions
\begin{equation}\label{eq:RLVE_bc}
\begin{array}{lclclcl}
h(t_0) &=& 79.248 \text{ km} &,& h(t_f) &=& 24.384 \text{ km}, \\[5pt]
\theta(t_0) &=& 0 \text{ deg} &,& \theta(t_f) &=& \text{Free}, \\[5pt]
\phi(t_0) &=& 0 \text{ deg} &,& \phi(t_f) &=& \text{Free}, \\[5pt]
v(t_0) &=& 7802.88 \text{ m/s} &,& v(t_f) &=& 762.0 \text{ m/s}, \\[5pt]
\gamma(t_0) &=& -1 \text{ deg} &,& \gamma(t_f) &=& -5 \text{ deg}, \\[5pt]
\psi(t_0) &=& 90 \text{ deg} &,& \psi(t_f) &=& \text{Free},
\end{array}
\end{equation}
where the terminal time, $t_f$, is free.
\section{Results}\label{section:results}
In this section, two studies are performed which analyze various forms of the hypersonic reusable launch vehicle entry (RLVE) optimal control problem presented in Section~\ref{section:probform}. The first study, given in Section~\ref{section:study_1}, serves as a verification and validation analysis of the SPOC method and contains two case scenarios. The two cases solve the RLVE optimal control problem with and without limits enforced on the control variables, while the single mixed state-control constraint (sensed acceleration) and both state-path constraints (heating rate and dynamic pressure) are enforced in both cases. The second study, given in Section~\ref{section:study_2}, examines the RLVE optimal control problem with Earth rotation included. The SPOC method is applied to obtain new results for this problem formulation, and is shown to detect and optimize the change in constraint structure due to the rotation of the Earth. Additionally, an analysis is performed on decreasing the maximum allowable stagnation point heating rate to further highlight that, without a priori knowledge of the solution structure, the SPOC method can detect and optimize changes in the structure of active state-path constraints. Finally, Section~\ref{section:control_discontinuity} includes a discussion on a feature of the SPOC method to capture small control discontinuities due to state-path constraints becoming active/inactive.

In the results to follow, the solutions obtained using the SPOC method are compared with a solution obtained using an $hp$-adaptive method which also employs Legendre-Gauss-Radau collocation (see Ref.~\cite{LiuRao2018} for further details). The solutions obtained using the mesh refinement method of Ref.~\cite{LiuRao2018} are henceforth referred to as $hp$-LGR. Additionally, to further verify and validate the results obtained using the SPOC method, comparisons will be made against the results provided in Ref.~\cite{MallTaheri2022} when possible, where it is noted Ref.~\cite{MallTaheri2022} utilizes an advanced indirect method to solve the same formulation of the RLVE problem. All results are obtained with the \textsf{MATLAB}$^\textsf{\textregistered}$ optimal control software $\mathbb{GPOPS-II}$\cite{PattersonRao2014} using the open source NLP solver IPOPT\cite{BieglerZavala2008} in full Newton mode with the linear solver MA57~\cite{ma57}, an NLP tolerance of $\epsilon_{NLP} = 1 \times 10^{-8}$, and a mesh refinement tolerance $\epsilon_{mesh} = 1 \times 10^{-7}$. The mesh refinement tolerance is the same for both the SPOC and $hp$-LGR solutions. All first and second derivatives were supplied to IPOPT using the built-in sparse central difference method in $\mathbb{GPOPS-II}$ developed in Ref.~\cite{PattersonRao2012}. Lastly, all computations were performed on a 2.9 GHz Intel Core i9 MacBook Pro running macOS Monterey Version 12.5 with 32 GB 2400 MHz DDR4 of RAM, using \textsf{MATLAB}$^\textsf{\textregistered}$ version R2023b (build 23.2.0.2365128).

\subsection{Study 1: Nonrotating Earth -- Comparison with Prior Results}\label{section:study_1}
In this study, two cases are considered to verify and validate the solutions obtained using the SPOC method. The first case considers the RLVE problem presented in Section~\ref{section:probform} \textit{without} any limits placed on the control variables. That is, the bounds in Eq.~\eqref{eq:controlconstraints} on the bank angle and angle of attack are set to $\sigma_{\min} = -\infty$ and $\alpha_{\max} = \infty$, respectively. The first case is included in this work as an analysis to verify the SPOC method can detect and optimize multiple state-path constraints active throughout the solution. Moreover, the results of the first case will be compared with the results obtained in the second study given in Section~\ref{section:study_2_comparison}.

The second case considers the RLVE problem presented in Section~\ref{section:probform} \textit{with} limits placed on the control variables. That is, the bounds in Eq.~\eqref{eq:controlconstraints} on the bank angle and angle of attack are set to $\sigma_{\min} = -75$ deg and $\alpha_{\max} = 19$ deg, respectively. The values of the control limits were chosen in order to directly compare results given in Ref.~\cite{MallTaheri2022}. The second case is included as an analysis to verify the SPOC method can also handle active control inequality path constraints while also detecting and optimizing any changes in the state-path constraint structure.

For both cases that follow, the initial mesh consists of $K = 30$ mesh intervals with $N_k = 5$, $k = 1,\ldots,K$, collocation points within each interval. Additionally, the initial guess for the initial mesh consists of a straight line connecting boundary conditions for states with boundary conditions specified at both the beginning and end of the original time domain, while a constant value is supplied for states with a single boundary condition specified at either end of the original time domain. As discussed in Sections~\ref{section:detection_decomposition} and~\ref{section:constraints_refinement}, the user-defined parameters required by the SPOC method are the detection tolerances $\epsilon_i$, $i = 1,\ldots,n_{sc}$, the detected switch time bound size parameters $\nu_i$, $i = 1,\ldots,n_{sc}$, and the maximum violation error tolerance, $\epsilon_c$. For both cases, the detection tolerances are set to $\epsilon_1 = 1 \times 10^{-5}$ and $\epsilon_2 = 1 \times 10^{-4}$, while the detected switch time bound size parameters are set to $\nu_1 = 0.5$, and $\nu_2 = 1$ corresponding to the heat rate and dynamic pressure constraints, respectively. It has been found that detection tolerances on the order of the square root of the NLP solver tolerance ($\sim \mathcal{O}(\sqrt{\epsilon_{NLP}})$) tend to return the best detection results \cite{ByczkowskiRao2023a}. Finally, the maximum relative constraint violation error tolerance is set equal to the mesh tolerance $\epsilon_c = \epsilon_{mesh} = 1 \times 10^{-7}$. 

\subsubsection{Case 1: Without Control Constraints}\label{section:case_1}
In this case, the inequality path constraints that are enforced on the system are the stagnation point heating rate, dynamic pressure, and sensed acceleration given by Eqs.~\eqref{eq:heatrate},~\eqref{eq:dynamicpressure}, and~\eqref{eq:sensedaccel}, respectively.
Table~\ref{table:Case1-detection_performance} provides the optimized entry and exit times for both state-path constraints computed upon each mesh iteration of the SPOC method. As the maximum mesh error is reduced, two observations can be made from Table~\ref{table:Case1-detection_performance}. The first observation is that the optimized entry and exit times for the active heating rate constraint converge to their respected values as the error is reduced. This observation is a result of the mesh being refined and, consequently, the search space for the detected times also being refined. Note, between each mesh iteration, the same mesh refinement technique being used to obtain the $hp$-LGR solution is being used within the SPOC method (see Refs.~\cite{LiuRao2018} and~\cite{ByczkowskiRao2024a} for further details). The second observation from Table~\ref{table:Case1-detection_performance} is that the active dynamic pressure constraint is not detected until the second mesh iteration and is subsequently optimized on the third mesh iteration. 
\begin{table}[b]
\centering
\caption{Optimized state-path constraint entry and exit times upon each mesh iteration of the SPOC method for Case 1 (without control constraints).}
\begin{tabular*}{\textwidth}{@{\extracolsep{\fill}}*{6}{c}}
\hline \hline
 \rule{0pt}{3ex}  Mesh iteration & \multicolumn{2}{c}{ Active arc 1, $\dot{Q} = \dot{Q}_{\max}$} & \multicolumn{2}{c}{Active arc 2, $q = q_{\max}$} & Max mesh error\\[3pt] \cline{1-1}  \cline{2-3} \cline{4-5} \cline{6-6}
 \rule{0pt}{3ex}  $M$ & Entry $[s]$ & Exit $[s]$ & Entry $[s]$ & Exit $[s]$ & $e_{\max}$  \\[3pt] \hline
 \rule{0pt}{3ex}0 &  -- & --  & -- & -- & 1.27$\times10^{-5}$ \\[3pt] 
 \rule{0pt}{3ex}1 & 236.16 & 694.81  & -- & -- & 9.52$\times10^{-6}$  \\[3pt] 
 \rule{0pt}{3ex}2 &  167.04 & 716.14  & -- & -- & 2.06$\times10^{-7}$\\[3pt] 
 \rule{0pt}{3ex}3 &  166.28 & 718.40  & 2085.45 & 2089.42 & 1.03$\times10^{-7}$ \\
 \rule{0pt}{3ex}4 &  165.73 & 716.50  & 2085.44 & 2089.32 & 9.54$\times10^{-8}$ \\[3pt] \hline \hline
\end{tabular*}
\label{table:Case1-detection_performance}
\end{table}
This observation emphasizes the importance of re-performing structure detection until all error tolerances have been met. To be clear, the detected times obtained using the structure detection method discussed in Section~\ref{section:detection_decomposition} are not shown in Table~\ref{table:Case1-detection_performance}, only the times that are optimized after solving the reformulated NLP are included. Also not provided in Table~\ref{table:Case1-detection_performance} is the observed maximum constraint violation tolerance on each iteration of the SPOC method, $\epsilon_c$, as this parameter was never violated throughout any of the mesh iterations. 

\begin{figure}[htb]
\begin{minipage}{.5\linewidth}
\centering
\subfloat[Dynamic pressure, $q(t)$, vs. time, $t$]{\label{fig:Case1_q}\includegraphics[scale=.41]{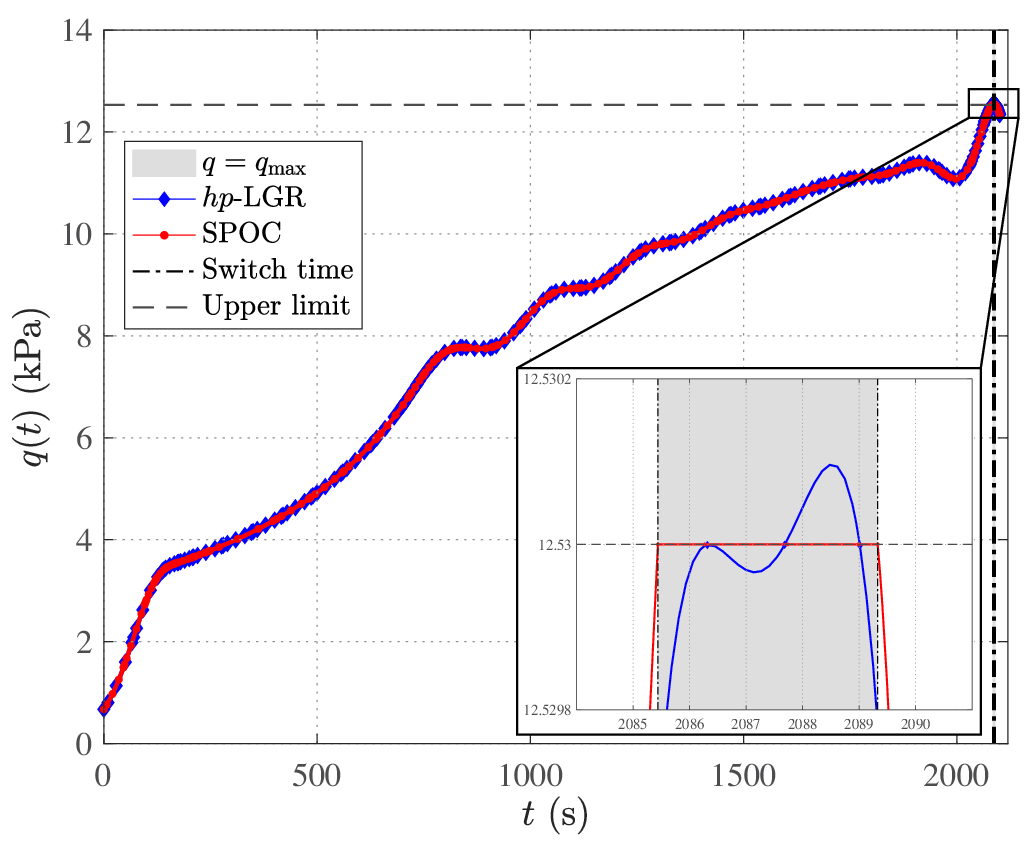}}
\end{minipage}%
\begin{minipage}{.5\linewidth}
\centering
\subfloat[Sensed acceleration load, $n(t)$, vs. time, $t$]{\label{fig:Case1_n}\includegraphics[scale=.41]{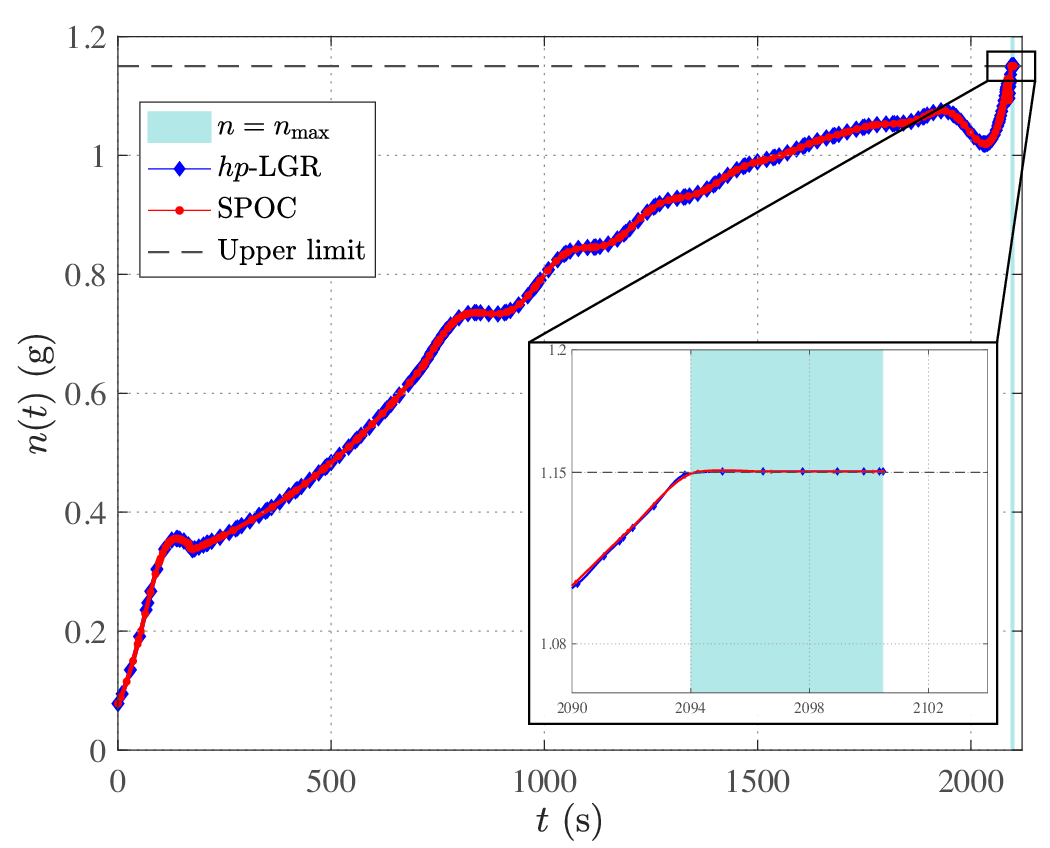}}
\end{minipage}\par\medskip
\centering
\subfloat[Stagnation point heating rate, $\dot{Q}(t)$, vs. time, $t$]{\label{fig:Case1_qdot}\includegraphics[scale=.41]{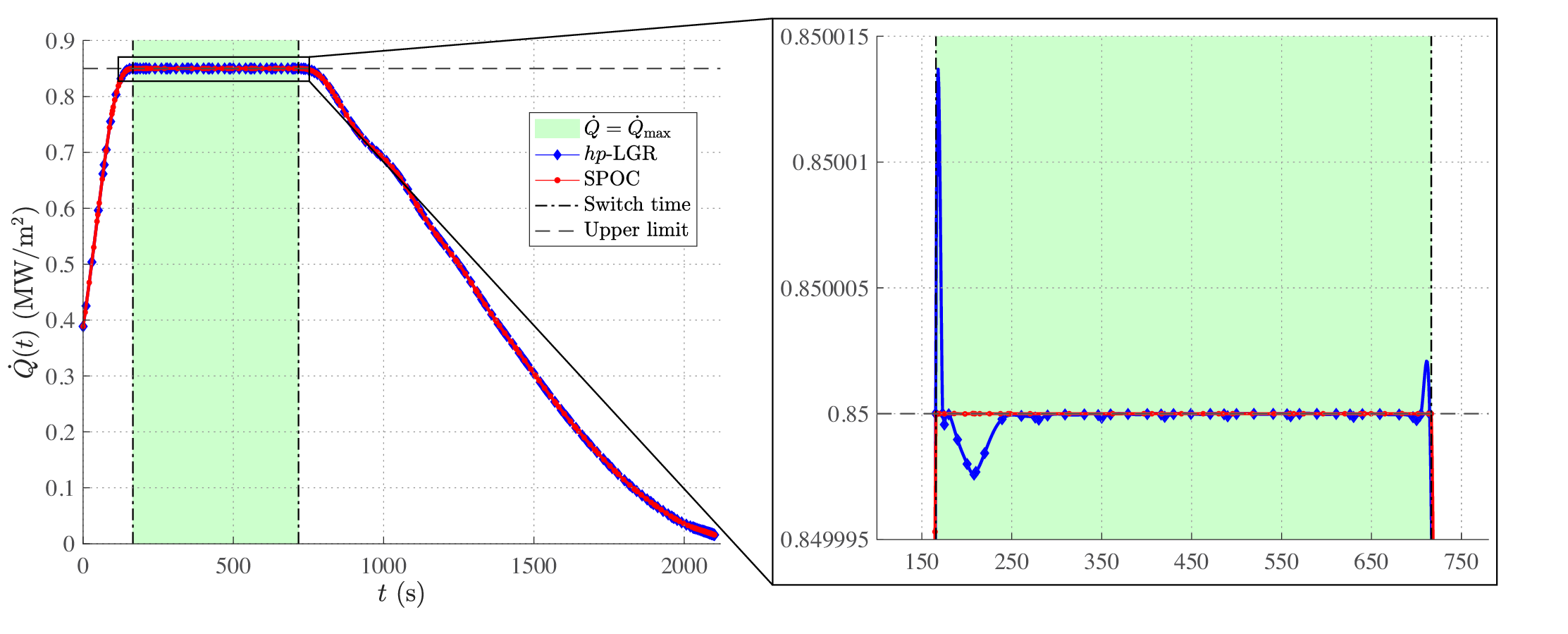}}
\caption{Optimal state-path constraint and mixed state-control constraint structures for Case 1 computed using the SPOC method compared against the $hp$-LGR solution.}
\label{fig:Case1_constraints}
\end{figure}

Next, Fig.~\ref{fig:Case1_constraints} provides the optimal dynamic pressure, sensed acceleration, and heating rate constraint profiles computed using the SPOC method alongside the $hp$-LGR solution. Specifically, the collocation points for both solutions are displayed along with an interpolated solution through the corresponding collocation points. A few observations can be made from Fig.~\ref{fig:Case1_constraints}. First, in Fig.~\ref{fig:Case1_q}, an enlarged view of the active dynamic pressure constraint shows that the interpolated SPOC solution is able to accurately capture the active arc, and unlike the interpolated $hp$-LGR solution, avoids any violations of the upper limit of $q_{\max} = 12.53 $kPa. 
A similar observation can be made from the enlarged view of the active heating rate arc provided in Fig.~\ref{fig:Case1_qdot}. Similar to previous findings (see Ref.~\cite{ByczkowskiRao2023b}), chattering behavior near the entry and exit times of the heat rate constraint arc in the $hp$-LGR solution is no longer present in the SPOC solution. Figure.~\ref{fig:Case1_n} shows that both the interpolated SPOC and $hp$-LGR solutions capture the active sensed acceleration arc (which is a mixed state-control path constraint and is not handled by the SPOC method as it is already a function of the angle of attack, $\alpha$). 

\begin{figure}[htb] 
\centering
\subfloat[Angle of attack, $\alpha(t)$, vs. time, $t$]{\label{fig:Case1_aoa}\includegraphics[scale=.4]{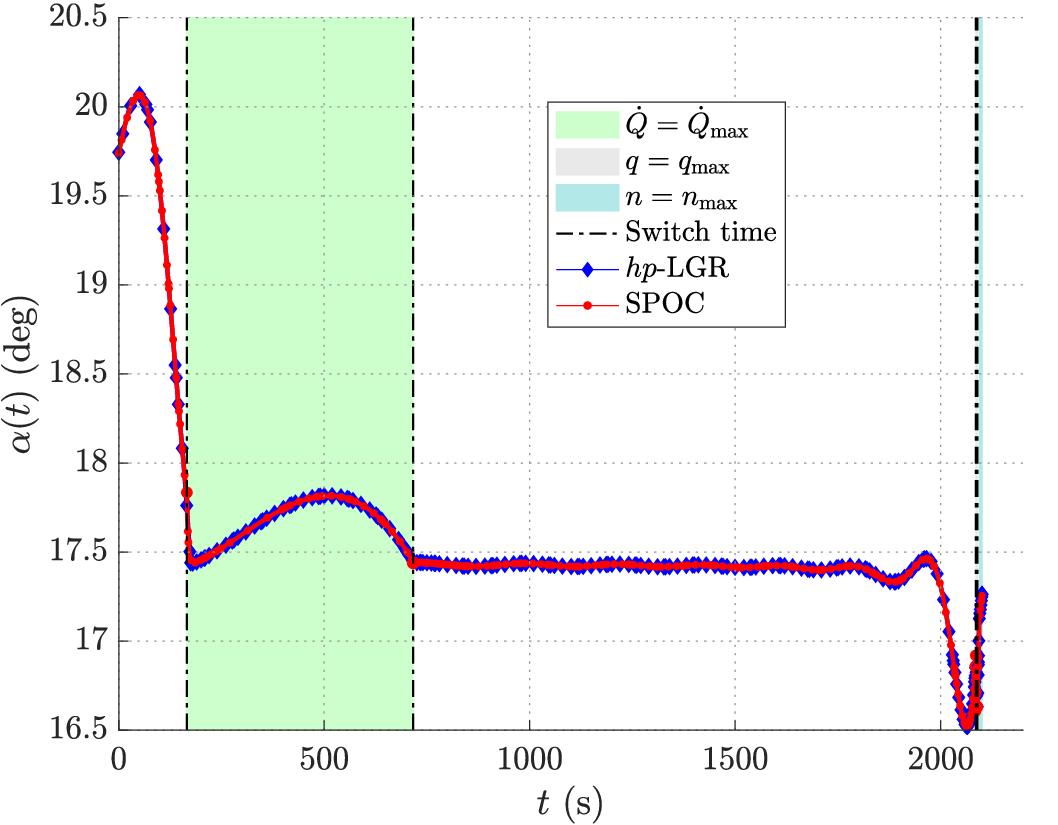}}
~~\subfloat[Bank angle, $\sigma(t)$, vs. time, $t$]{\label{fig:Case1_bank}\includegraphics[scale=.4]{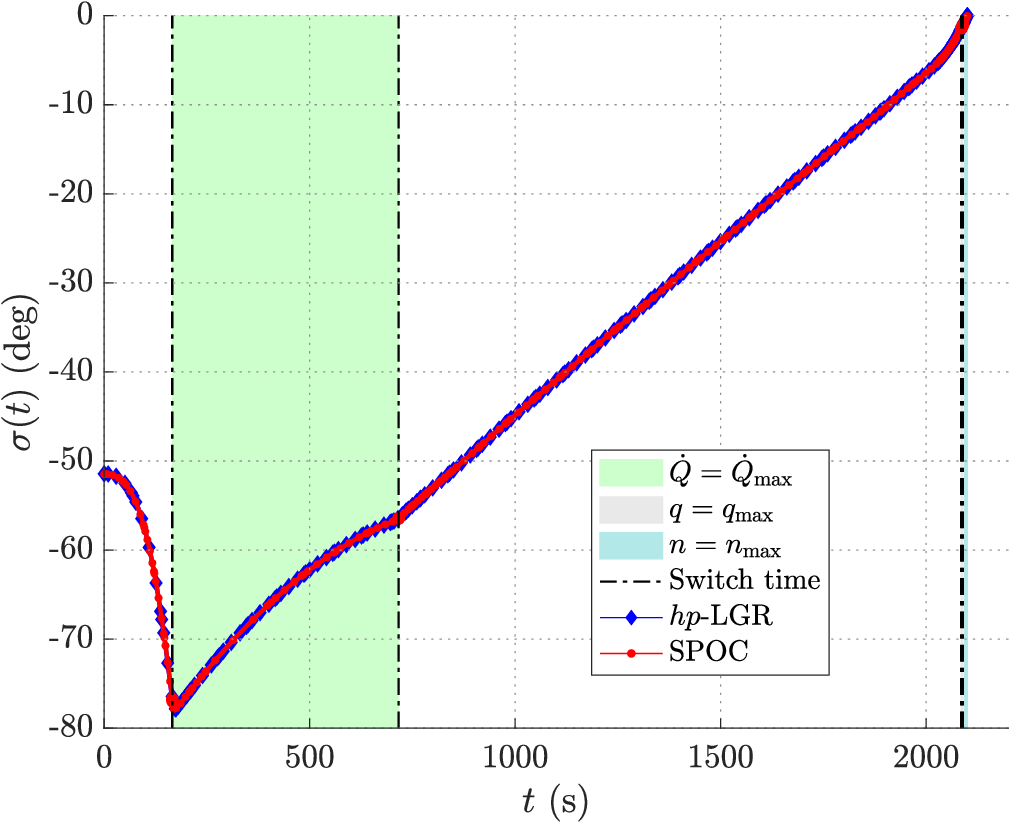}}
\caption{Optimal angle of attack, $\alpha(t)$, and bank angle, $\sigma(t)$, profiles for Case 1 using the SPOC method compared against an $hp$-LGR solution.\label{fig:Case1_control}}
\end{figure}

Lastly, for Case 1, Fig.~\ref{fig:Case1_control} provides the optimal control components computed using the SPOC method compared against the $hp$-LGR solution. Rapid changes in both control variables are observed as the heating rate constraint becomes active. Additionally, the angle of attack experiences another rapid change near the activation of the dynamic pressure constraint. The vertical dashed lines appearing in both Figs.~\ref{fig:Case1_constraints} and~\ref{fig:Case1_control} correspond to the optimized entry and exit times provided in the final row of Table~\ref{table:Case1-detection_performance}.

\subsubsection{Case 2: With Control Constraints}\label{section:case_2}
In this case, the control inequality path constraints on the angle of attack and bank angle given by Eq.~\eqref{eq:controlconstraints} are enforced in addition to the inequality path constraints enforced in Case 1. The limits on all inequality path constraints are chosen to match those used in Ref.~\cite{MallTaheri2022}, which were chosen such that every constraint becomes active throughout the solution.
Table~\ref{table:Case2-detection_performance} provides the optimized entry and exit times for both state-path constraints obtained upon each mesh iteration of the SPOC method. As the maximum mesh error is reduced, a few observations can be made from Table~\ref{table:Case2-detection_performance}. The first observation is there now exists two separate active heating rate arcs. Similar to the dynamic pressure constraint, the first active heating rate constraint is not detected on the initial mesh. Specifically, on the initial mesh ($M = 0$), the first heating rate constraint only appears as a single touch point, but as the maximum mesh error is reduced, the first heating rate constraint becomes active over a non-zero interval of time. 
\begin{table}[htb]
\centering
\caption{Optimized state-path constraint entry and exit times using the SPOC method for Case 2.}
\begin{tabular*}{\textwidth}{@{\extracolsep{\fill}}*{8}{c}}
\hline \hline
 \rule{0pt}{3ex}   & \multicolumn{2}{c}{Active arc 1, $\dot{Q} = \dot{Q}_{\max}$} & \multicolumn{2}{c}{ Active arc 2, $\dot{Q} = \dot{Q}_{\max}$} & \multicolumn{2}{c}{Active arc 3, $q = q_{\max}$} & Max mesh error \\[3pt] \cline{2-3} \cline{4-5} \cline{6-7} \cline{8-8} 
 \rule{0pt}{3ex} $M$ & Entry $(s)$ & Exit $(s)$ & Entry $(s)$ & Exit $(s)$ & Entry $(s)$ & Exit $(s)$ & $e_{\max}$ \\[3pt] \hline
 \rule{0pt}{3ex}0 & -- & -- &  -- & --  & -- & -- & 1.27$\times10^{-5}$ \\[3pt]  
 \rule{0pt}{3ex}1 & -- & -- &  416.34 & 724.47  & -- & -- & 1.08$\times10^{-5}$ \\[3pt] 
 \rule{0pt}{3ex}2 & 168.19 & 169.66 &  411.71 & 728.88  & -- & -- & 1.51$\times10^{-6}$\\[3pt] 
 \rule{0pt}{3ex}3 & 167.05 & 169.57 &  412.13 & 728.45  & 2095.42 & 2096.07 & 2.40$\times10^{-7}$\\[3pt] 
 \rule{0pt}{3ex}4 & 167.03 & 169.70 &  411.16 & 728.95  & 2095.41 & 2099.01 & 7.95$\times10^{-8}$\\[3pt] \hline \hline
\end{tabular*}
\label{table:Case2-detection_performance}
\end{table}
\begin{figure}[htb!]
\centering
\subfloat[Enlarged view of heating rate profile for Case 2]{\label{fig:Case2_qdot}\includegraphics[scale=.41]{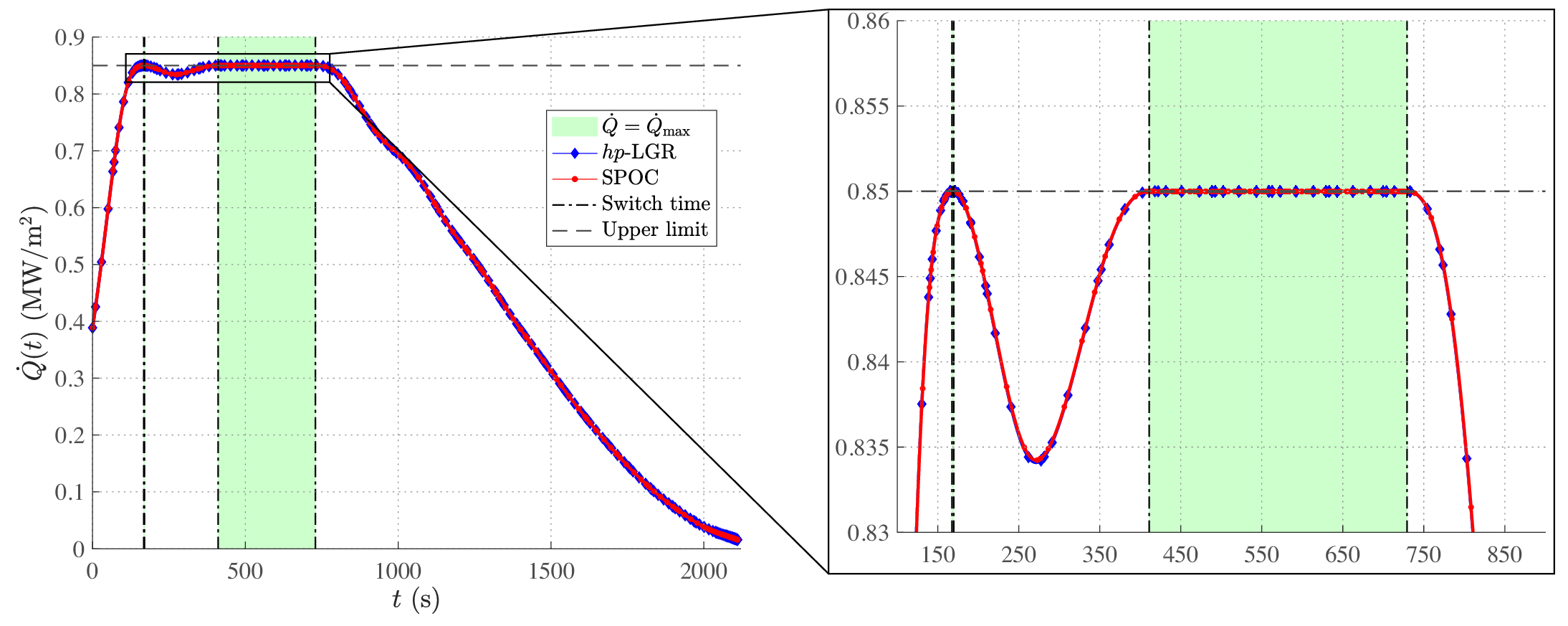}}\par\medskip
\begin{minipage}{.5\linewidth}
\centering
\subfloat[Enlarged view of first active heating rate arc]{\label{fig:Case2_Qdot1}\includegraphics[scale=.41]{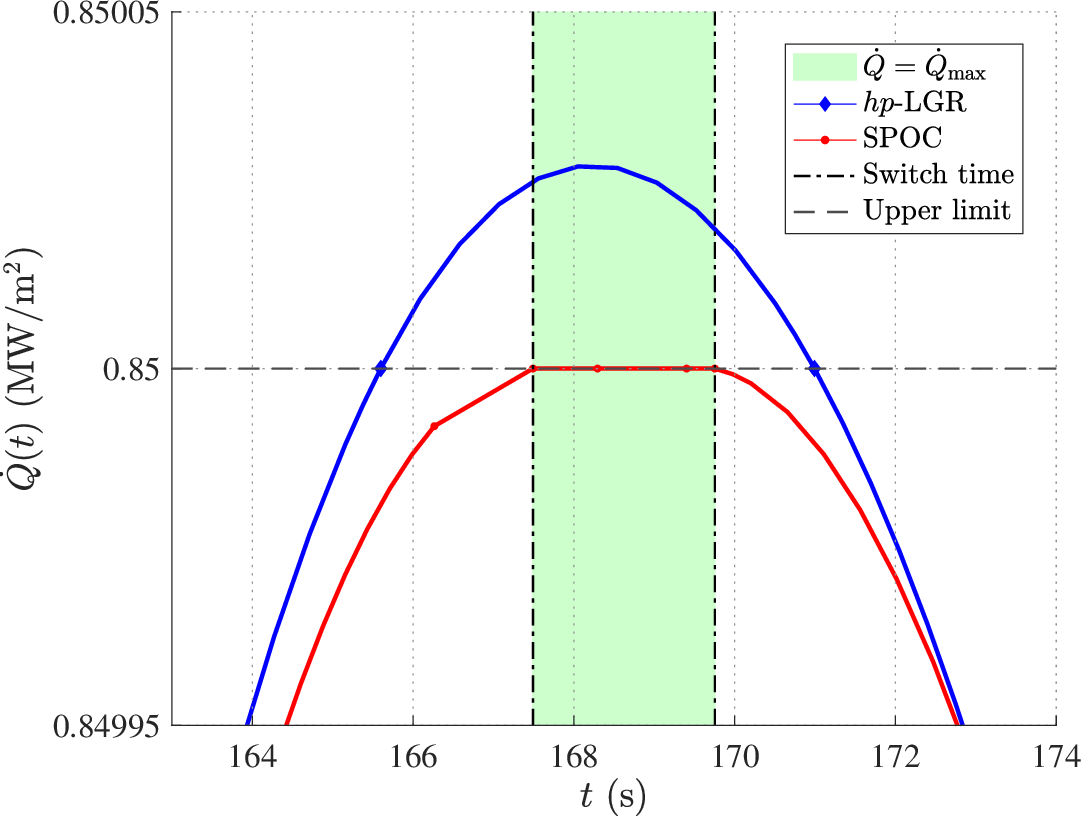}}
\end{minipage}%
\begin{minipage}{.5\linewidth}
\centering
\subfloat[Enlarged view of second active heating rate arc]{\label{fig:Case2_Qdot2}\includegraphics[scale=.41]{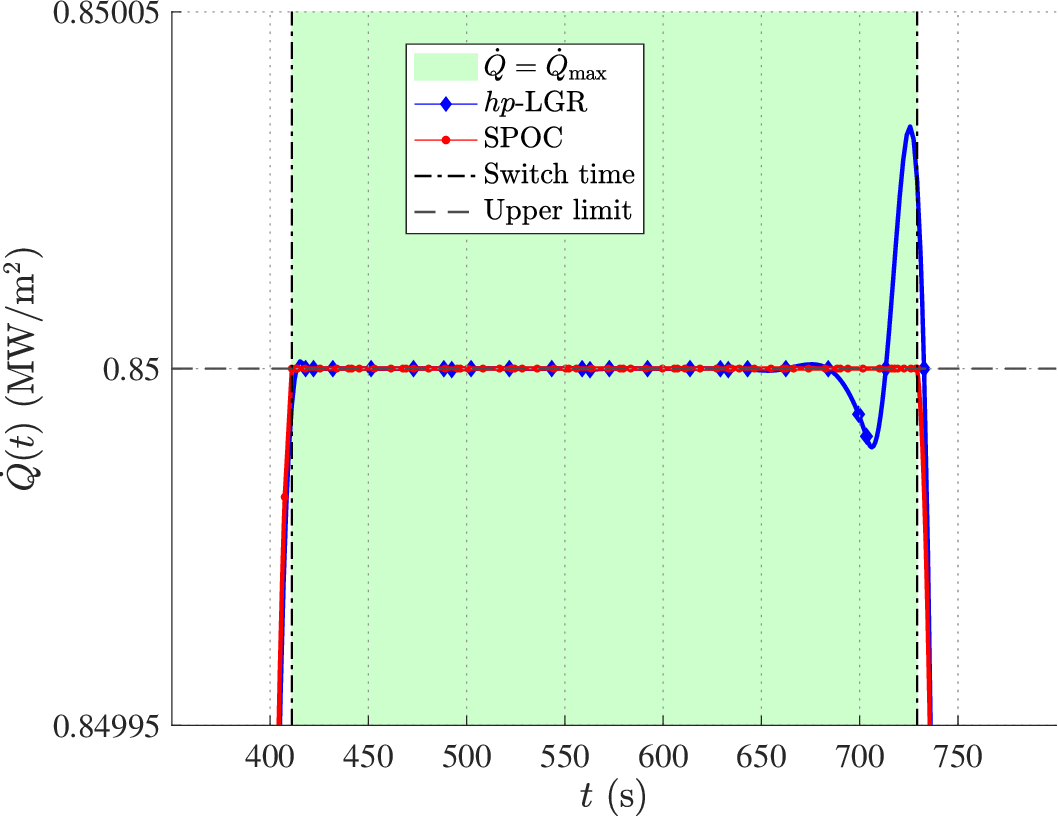}}
\end{minipage}
\caption{Optimal heating rate constraint structure for Case 2 computed using the SPOC method compared against the $hp$-LGR solution.}
\label{fig:Case2_heatrate}
\end{figure}
The next observation (which is similar to that made in Case 1) is as the maximum mesh error is reduced, the detected entry and exit times converge to their respected values. The convergence behavior of the detected/optimized switch times is due to the mesh being refined and the optimal constraint structure being accurately captured. The final observation, which again highlights the importance of re-performing structure detection, is that the active dynamic pressure constraint arc is not detected until the second to last mesh. The times listed in Tables~\ref{table:Case1-detection_performance} and~\ref{table:Case2-detection_performance} are the optimized times after solving the reformulated NLP, meaning the entry and exit times are detected on the previous mesh and optimized when solving the reformulated problem. Finally, note the SPOC method only performs structure detection of the heat rate and dynamic pressure constraints. Structure detection is only performed on these two constraints because they are purely functions of the state-variables.

\begin{figure}[htb]
\begin{minipage}{.5\linewidth}
\centering
\subfloat[Dynamic pressure, $q(t)$, vs. time, $t$]{\label{fig:Case2_q}\includegraphics[scale=.41]{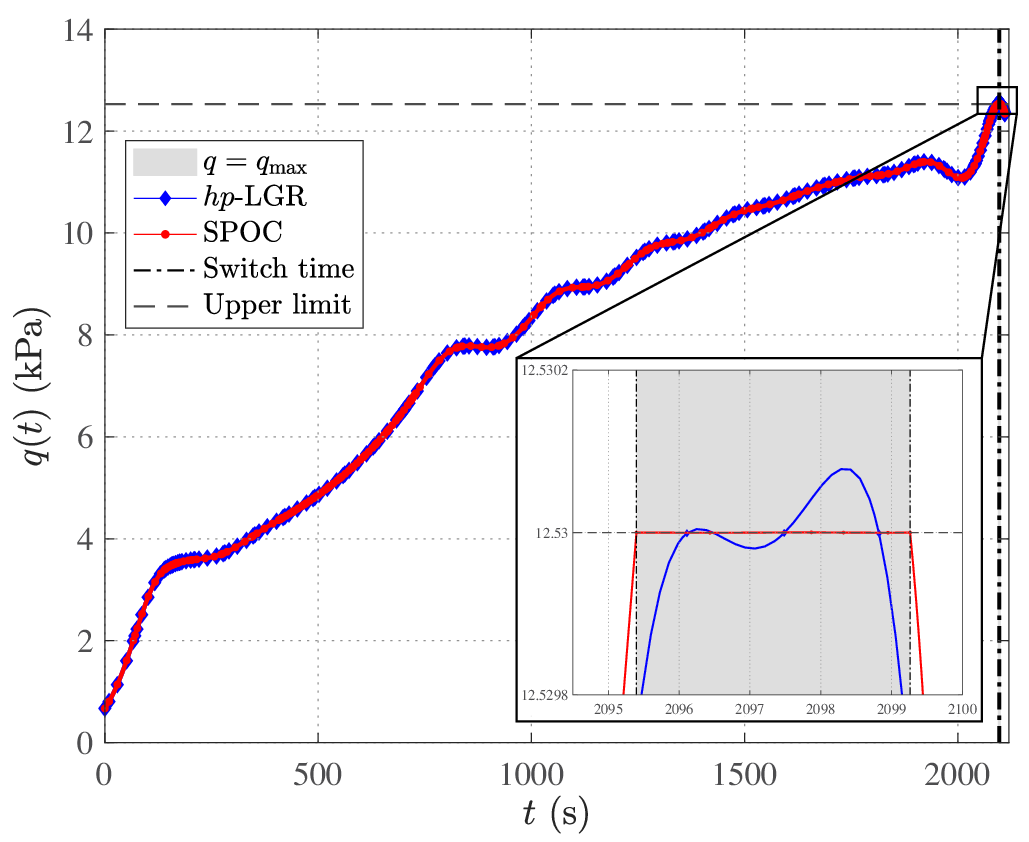}}
\end{minipage}%
\begin{minipage}{.5\linewidth}
\centering
\subfloat[Sensed acceleration load, $n(t)$, vs. time, $t$]{\label{fig:Case2_n}\includegraphics[scale=.41]{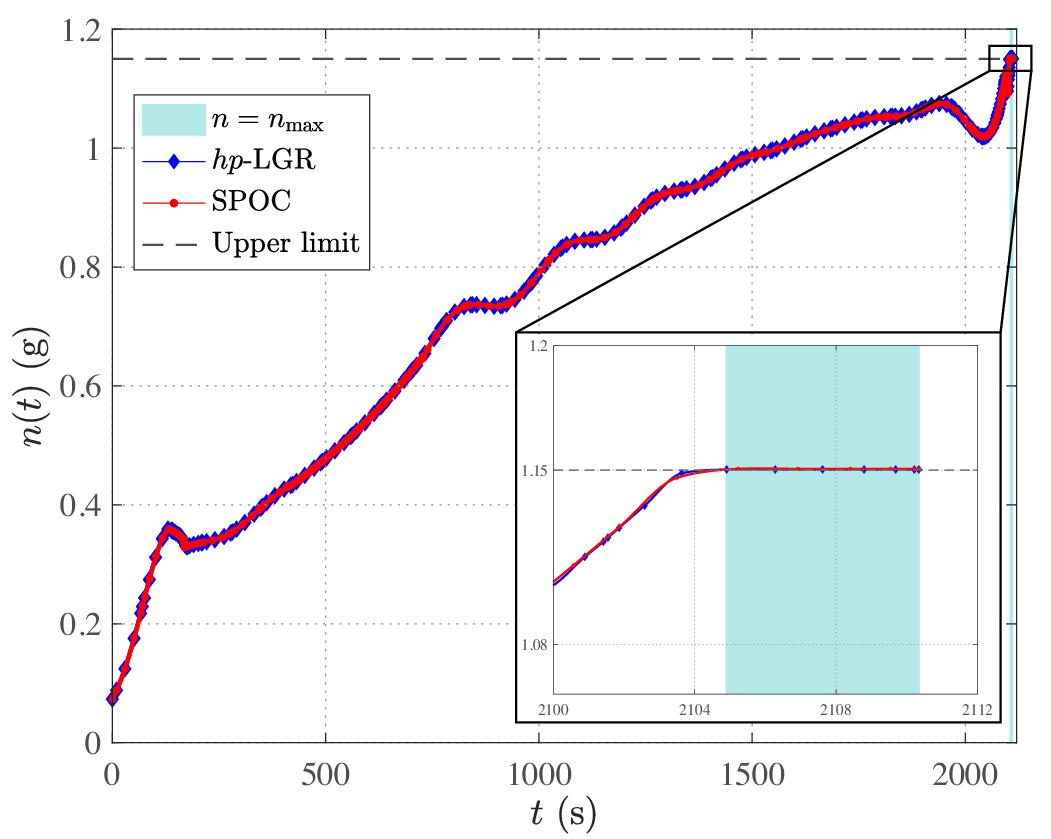}}
\end{minipage}\par\medskip
\caption{Optimal dynamic pressure and sensed acceleration constraint structures for Case 2 computed using the SPOC method compared against the $hp$-LGR solution.}
\label{fig:Case2_constraints}
\end{figure}

Next, Figs.~\ref{fig:Case2_heatrate} and~\ref{fig:Case2_constraints} provide the optimal heating rate, dynamic pressure, and sensed acceleration constraint profiles computed using the SPOC method alongside the $hp$-LGR solution. Figure~\ref{fig:Case2_qdot} shows the separate active heat rate arcs as a result of enforcing the control constraints. Comparing the interpolated SPOC solution to the interpolated $hp$-LGR solution shown in Figs.~\ref{fig:Case2_Qdot1} and~\ref{fig:Case2_Qdot2}, it is found that the SPOC method is able to remove any violations of the upper limit and reduce chattering near the entry and exit times. Similar to Case 1, Fig.~\ref{fig:Case2_q} provides an enlarged view of the active dynamic pressure constraint which shows that the interpolated SPOC solution is able to accurately capture the active arc, while the interpolated $hp$-LGR solution still violates the upper limit of $q_{\max} = 12.53 $kPa. Figure~\ref{fig:Case2_n} shows that both the interpolated SPOC and $hp$-LGR solutions capture the active sensed acceleration arc. Recall, the vertical dashed lines in Figs.~\ref{fig:Case2_heatrate} and~\ref{fig:Case2_q} correspond to the optimized times provided in the last row of Table~\ref{table:Case2-detection_performance}. 

\begin{figure}[htb] 
\centering
\subfloat[Angle of attack, $\alpha(t)$, vs. time, $t$]{\label{fig:Case2_aoa}\includegraphics[scale=.4]{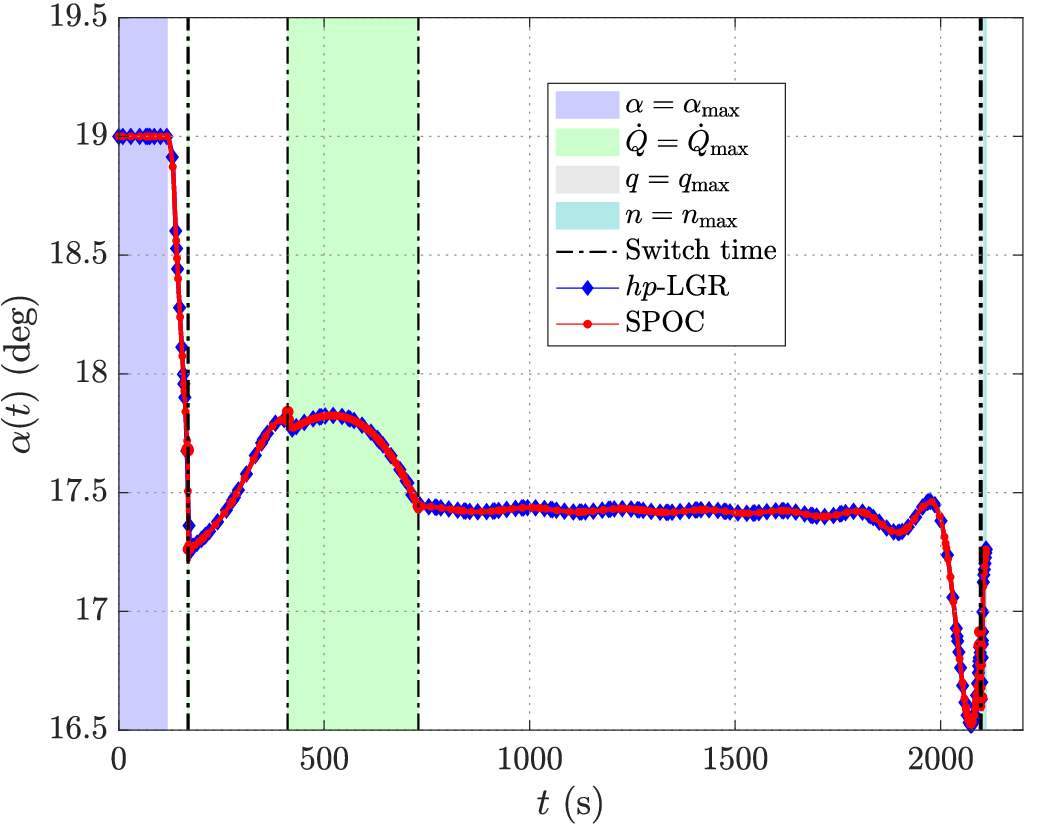}}
~~\subfloat[Bank angle, $\sigma(t)$, vs. time, $t$]{\label{fig:Case2_bank}\includegraphics[scale=.4]{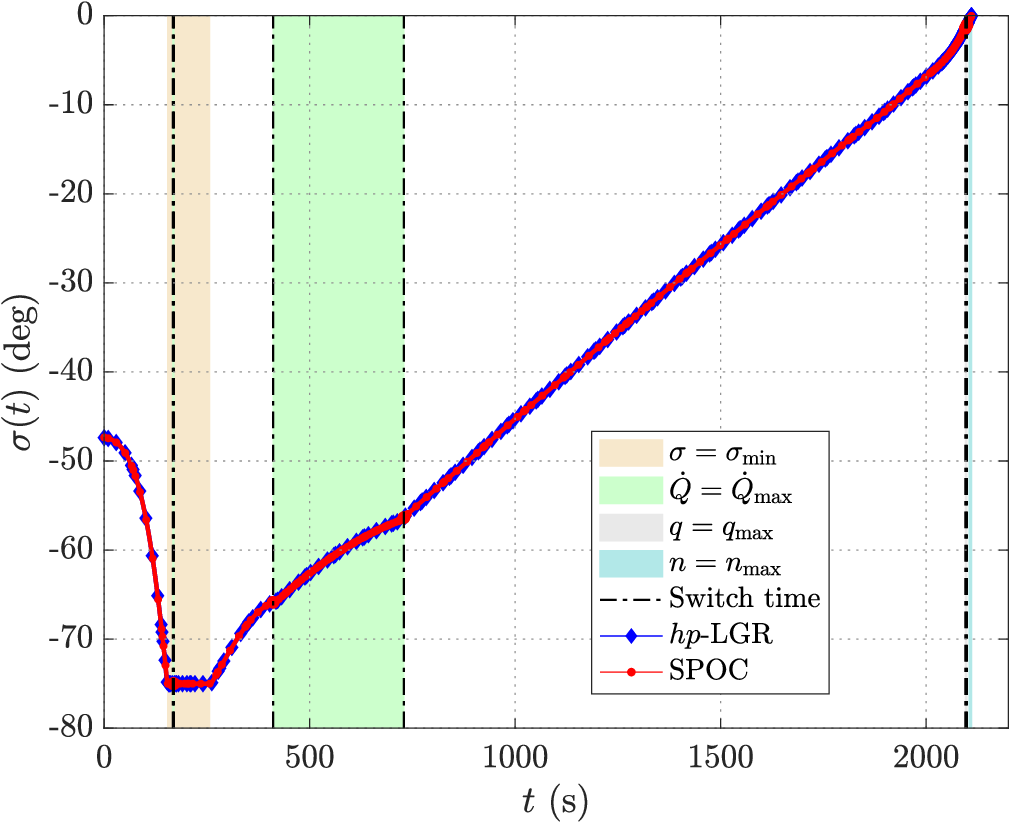}}
\caption{Optimal angle of attack, $\alpha(t)$, and bank angle, $\sigma(t)$, profiles for Case 2 using the SPOC method compared against an $hp$-LGR solution.\label{fig:Case2_control}}
\end{figure}

Lastly for Case 2, Fig.~\ref{fig:Case2_control} provides a comparison of the control components obtained using the SPOC method against the $hp$-LGR solution. It is found that the optimal angle of attack starts at the upper limit of $\alpha_{\max} = 19$ deg while the bank angle constraint becomes active thereafter. In particular, the bank angle constraint becomes active at a time that is close in proximity to the first heating rate constraint arc. 

\subsubsection{Case Comparison and Discussion}
Table~\ref{table:Case_comparison} provides the entry and exit times of each state-path constraint arc along with the optimal objective value (final crossrange) for both cases. Comparisons are made against the $hp$-LGR solution while an additional comparison is made against the results provided in Ref.~\cite{MallTaheri2022} for Case 2. It can be found that the SPOC solution is in close agreement with the solutions obtained by the  $hp$-adaptive collocation method as well as the advance indirect regularization method used in Ref.~\cite{MallTaheri2022}, further validating the solution obtained using the SPOC method. Specifically, for Case 2, the optimized entry and exit times of each state-path constraint are in close agreement with those provided in Ref.~\cite{MallTaheri2022}. Note, the inclusion of the control inequality path constraints does not significantly impact the final crossrange achieved, though it provided a good analysis to compare the change in state-path constraint structures. That is, the study of both cases was performed to highlight the SPOC method automatically detecting and optimizing the change in constraint structure. \par 
\begin{table}[htb]
\centering
\caption{Comparison of entry and exit times, and optimal objective value for both cases.}
\begin{tabular*}{\textwidth}{@{\extracolsep{\fill}}*{9}{c}}
\hline \hline
 \rule{0pt}{3ex}  & & \multicolumn{2}{c}{Arc 1, $\dot{Q} = \dot{Q}_{\max}$} & \multicolumn{2}{c}{ Arc 2, $\dot{Q} = \dot{Q}_{\max}$} & \multicolumn{2}{c}{Arc 3, $q = q_{\max}$} & Objective \\[3pt] \cline{3-4}\cline{5-6}\cline{7-8}\cline{9-9} 
 \rule{0pt}{3ex} Case & Method & Entry $(s)$ & Exit $(s)$ & Entry $(s)$ & Exit $(s)$ & Entry $(s)$ & Exit $(s)$ & $\phi(t_f)$ (deg) \\[3pt] \hline
 \multirow{2}{*}{1} \rule{0pt}{3ex} & $hp$-LGR & -- &  -- &  165.35 & 714.74  & 2086.32 & 2089.02 & 33.99 \\
 \rule{0pt}{3ex}  & SPOC & -- &  -- &  165.73 & 716.50  & 2085.44 & 2089.32 & 33.99 \\[3pt] 
\hline
\multirow{4}{*}{2} \rule{0pt}{3.1ex} & $hp$-LGR & 165.60 & 170.99 &  413.84 & 732.74  & 2096.11 & 2098.82 & 33.99 \\
                            \rule{0pt}{3.1ex}  & SPOC & 167.03 & 167.70 &  411.16 & 728.95  & 2095.41 & 2099.01 & 33.99  \\
                            \rule{0pt}{3.1ex} & Ref.~\cite{MallTaheri2022} & 167.15 &  168.83 &  417.40 & 724.02  & 2095.64 & 2098.12 & 33.99 \\[3pt]
\hline \hline
\end{tabular*}
\label{table:Case_comparison}
\end{table}
\begin{figure}[htb!] 
\centering
\subfloat[Energy plot comparing the unconstrained and fully constrained (Case 2) solutions]{\label{fig:energyplot}\includegraphics[scale=.4]{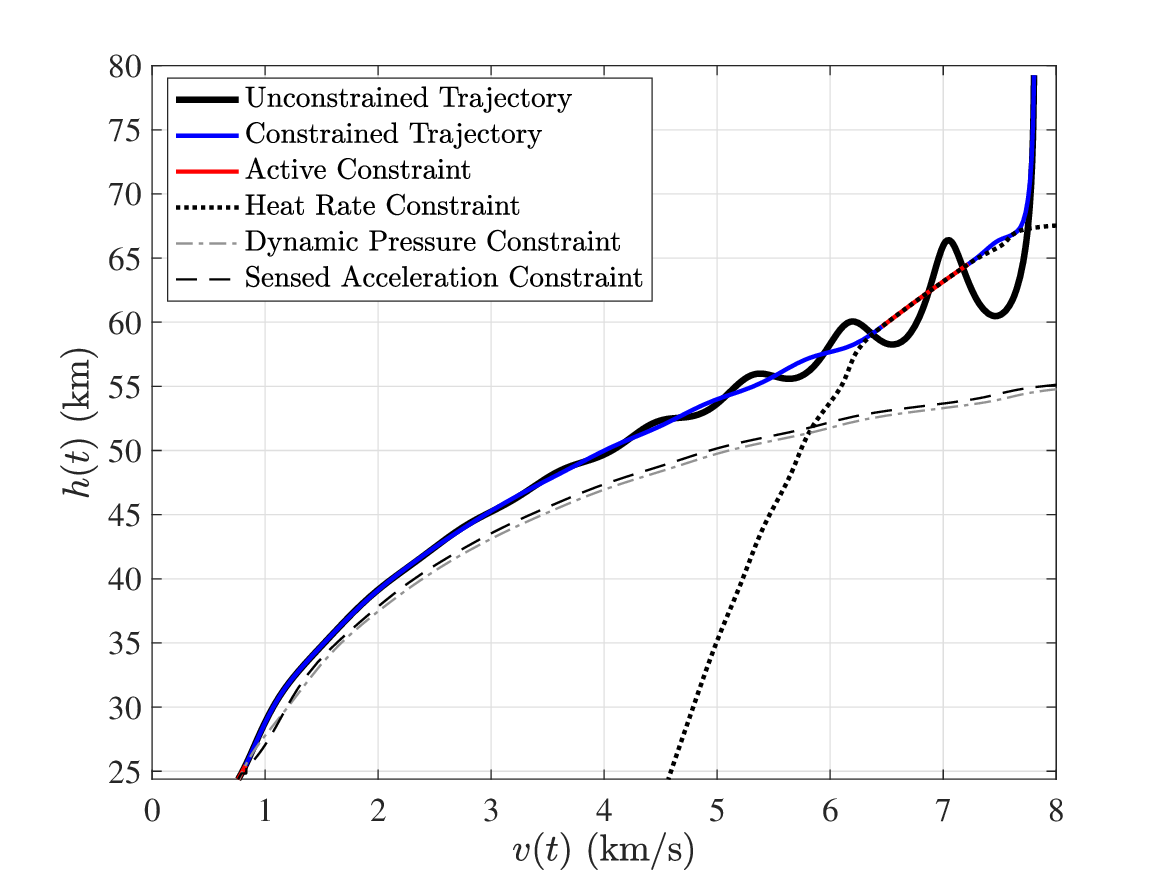}}
~~\subfloat[Three-dimensional entry trajectories]{\label{fig:3Dtraj}\includegraphics[scale=.4]{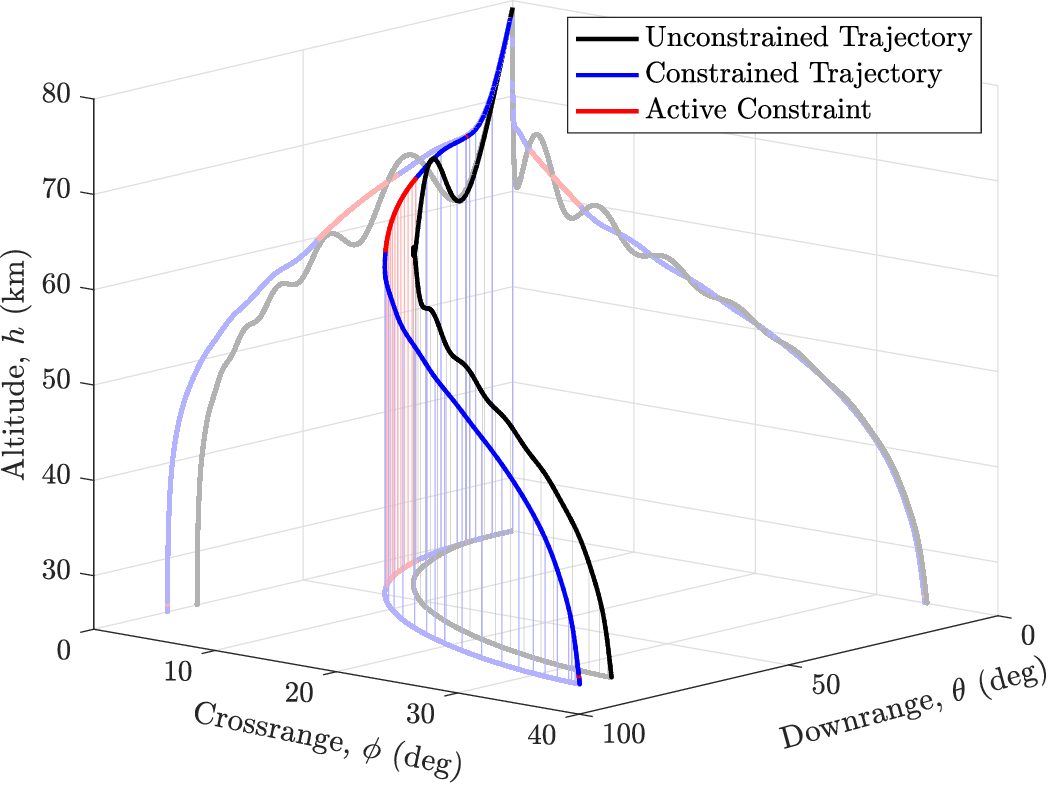}}
\caption{Comparison between the unconstrained and fully constrained (Case 2) entry trajectories.\label{fig:comparison}}
\end{figure}
The final values of the downrange and time of flight for Case 2 obtained using the SPOC method are 82.41 deg and 2110.37 s, respectively, while the corresponding results given in Ref~\cite{MallTaheri2022} are 82.42 deg and 2110.50 s, respectively, and that on the $hp$-LGR solution are 82.41 deg and 2110.34 s, respectively. Next, Fig.~\ref{fig:comparison} provides a comparison between the unconstrained solution and the fully constrained (Case 2) solution. Figure~\ref{fig:energyplot} shows the energy plot for both the unconstrained and constrained solution, where it can be seen that the heating rate constraint becomes active in a region where the vehicle would typically exhibit large dives into and out of the atmosphere (which is a phenomena also known as phugoid oscillations). The impact of the heating rate constraint smooths out the phugoid oscillations and the vehicle flies in a near equilibrium glide configuration throughout the upper part of the atmosphere. Also shown in Fig.~\ref{fig:energyplot} is the dynamic pressure constraint and sensed acceleration constraint becoming active towards the end of the trajectory. Figure~\ref{fig:3Dtraj} provides the three-dimensional entry trajectory for the unconstrained and fully constrained RLVE problem formulations. In addition to smoothing out the phugoid oscillations present in the unconstrained trajectory, it is found the inclusion of constraints result in the vehicle achieving a larger final downrange (longitude) distance while achieving a similar final crossrange (latitude) distance.

\subsection{Study 2: Rotating Earth -- New Results}\label{section:study_2}
In this study, a new problem formulation of the RLVE problem is solved with the SPOC method. Specifically, the rotation of the Earth is included in the system model, and the SPOC method is used to generate new results by detecting and optimizing any changes in the active state-path constraint structure. The remainder of this section is organized as follows. First, the equations of motion presented in Section~\ref{section:probform} are adjusted to include the rotation of the Earth. Next, Section~\ref{section:study_2_comparison} includes results obtained with Earth rotation and draws a comparison with the results presented in Section~\ref{section:case_1}. Lastly, Section~\ref{section:study_2_variation} performs an analysis on varying the maximum allowable stagnation point heating rate to study the change in the structure of the solution. The SPOC method is used to algorithmically detect and optimize the change in the solution structure without any prior knowledge of the number and location of active state-path constraints.  In the results to follow, the initial mesh consisted of $K = 10$ mesh intervals with $N_k = 4$ collocation points within each mesh interval. The solution obtained in Case 1 is supplied as the initial guess and the variables used within the SPOC method are the same values used to obtain the results presented in Section~\ref{section:study_1}.

\subsubsection{Comparison with Case 1}\label{section:study_2_comparison} 
To account for the rotation of the Earth, the kinetic equations given in Eqs.~\eqref{eq:vdot}-\eqref{eq:psidot} become \cite{VinhCulp1980}
\begin{subequations}\label{eq:EOM_Earthrotation}
\begin{align}
\dot{v} &= -D -g\sin \gamma + r \omega_e^2 \cos \phi (\sin \gamma \cos \phi - \cos \gamma \sin \phi \cos \psi ), \label{eq:vdotER}\\[5pt]
\dot{\gamma} &= \dfrac{L \cos \sigma}{v} + \cos \gamma \left( \dfrac{v}{r} - \dfrac{g}{v}  \right) + 2 \omega_e \cos \phi \sin \psi + \dfrac{ r \omega_e^2}{v} \cos \phi (\cos \gamma \cos \phi + \sin \gamma \sin \phi \cos \psi), \label{eq:gammadotER}\\[5pt]
\dot{\psi} &=  \dfrac{L \sin \sigma}{v \cos \gamma} + \dfrac{v}{r} \cos \gamma \sin \psi \tan \phi - 2 \omega_e (\tan \gamma \cos \phi \cos \psi - \sin \phi) + \dfrac{r \omega_e^2}{v \cos \gamma} \sin \phi \cos \phi \sin \psi ,\label{eq:psidotER}
\end{align}
\end{subequations}
where $\omega_e$ is the rotation rate of the Earth, $\theta$ is now the Earth-relative longitude, $\phi$ is the geocentric latitude, $v$ is now the Earth-relative speed, $\gamma$ is now the Earth-relative flight path angle, and $\psi$ is now the Earth-relative azimuth angle (where Ref.~\cite{VinhCulp1980} uses heading angle). Note, the kinematic equations given by Eqs.~\eqref{eq:rdot}-\eqref{eq:phidot} remain the same across problem formulations.

\begin{figure}[htb]
\begin{minipage}{.5\linewidth}
\centering
\subfloat[Dynamic pressure, $q(t)$, vs. time, $t$]{\label{fig:study2_q}\includegraphics[scale=.40]{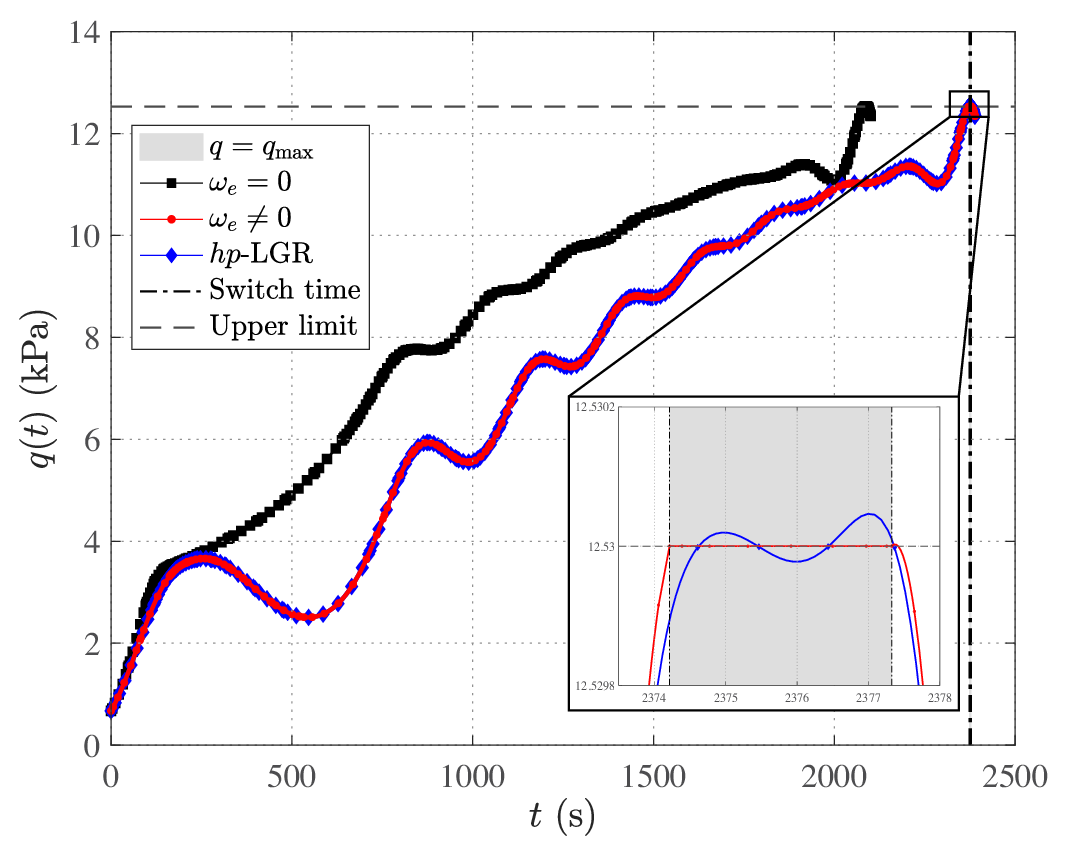}}
\end{minipage}%
\begin{minipage}{.5\linewidth}
\centering
\subfloat[Sensed acceleration load, $n(t)$, vs. time, $t$]{\label{fig:study2_n}\includegraphics[scale=.40]{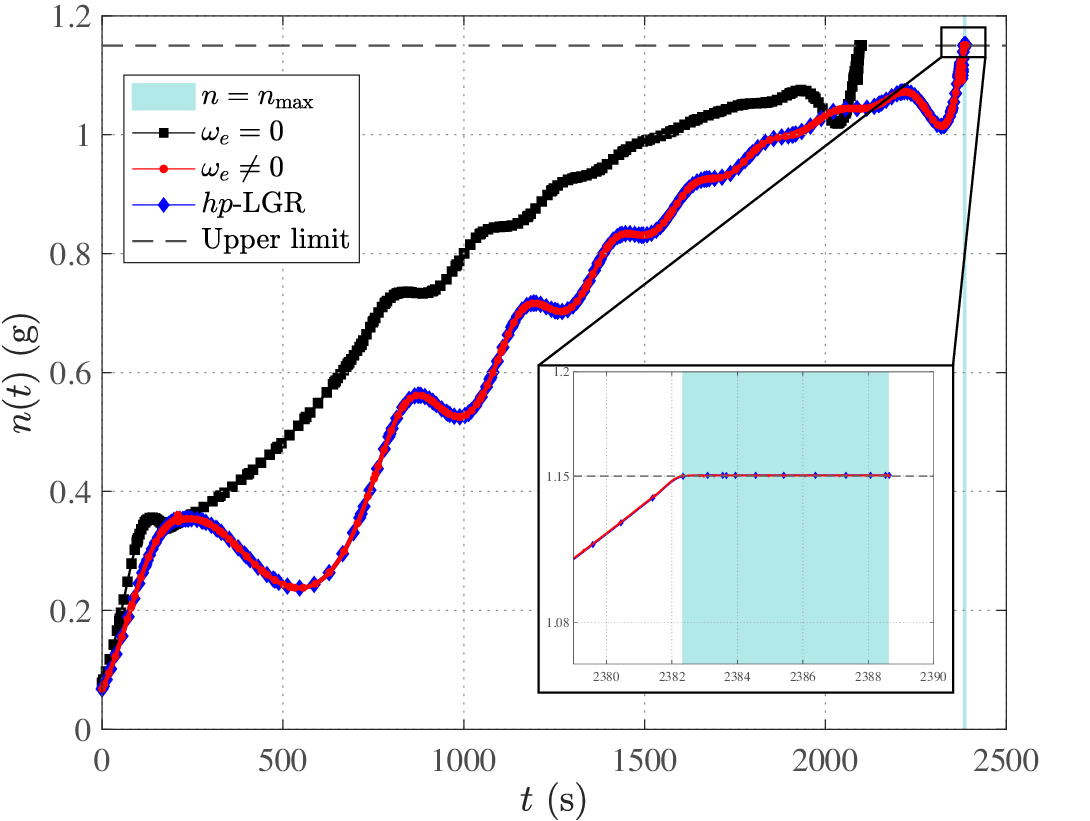}}
\end{minipage}\par\medskip
\begin{minipage}{.50\linewidth}
\centering
\subfloat[Stagnation point heating rate, $\dot{Q}(t)$, vs. time, $t$]{\label{fig:study2_qdot}\includegraphics[scale=.40]{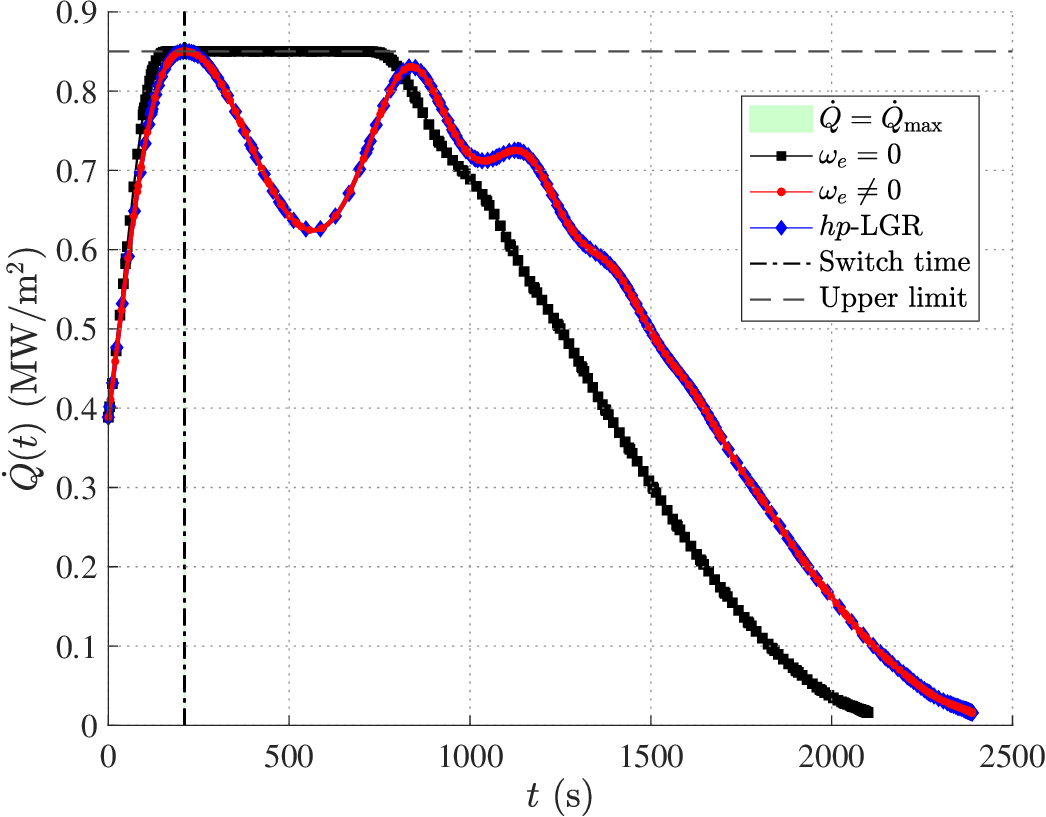}}
\end{minipage}
\begin{minipage}{.45\linewidth}
\centering
\subfloat[Enlarged view of active heating rate arc]{\label{fig:study2_qdotzoom}\includegraphics[scale=.38]{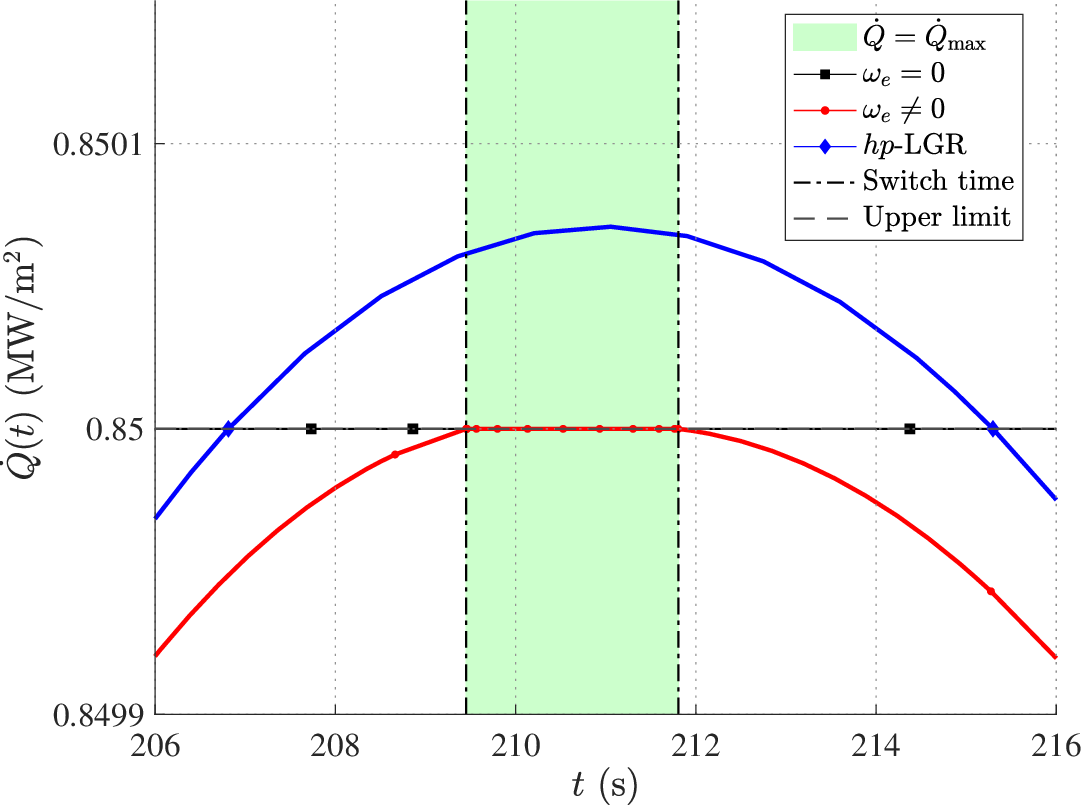}}
\end{minipage}
\caption{Optimal state-path constraint and mixed state-control constraint profiles for Study 2 compared against the results obtained in Case 1.}
\label{fig:Study2_constraints}
\end{figure}

Figure~\ref{fig:Study2_constraints} includes the optimal state-path constraint and mixed state-control constraint structures for the RLVE vehicle with and without Earth rotation. 
It is seen in Figs.~\ref{fig:study2_qdot} and~\ref{fig:study2_qdotzoom} that, with the inclusion of the rotation of the Earth, the duration of the active heating rate constraint is significantly reduced, but not completely removed, as the SPOC method detects an active arc of approximately two seconds long. Figures~\ref{fig:study2_q} and~\ref{fig:study2_n} show that the dynamic pressure and sensed acceleration constraints follow a similar trend to the non-rotating Earth model. Lastly, it is shown in Fig.~\ref{fig:Study2_constraints} that the SPOC solution is in close agreement with the solution obtained using an $hp$-adaptive method, while improving the accuracy of the computed active state-path constraints.

\begin{figure}[htb!] 
\centering
\subfloat[Angle of attack, $\alpha(t)$, vs. time, $t$]{\label{fig:Study2_aoa}\includegraphics[scale=.4]{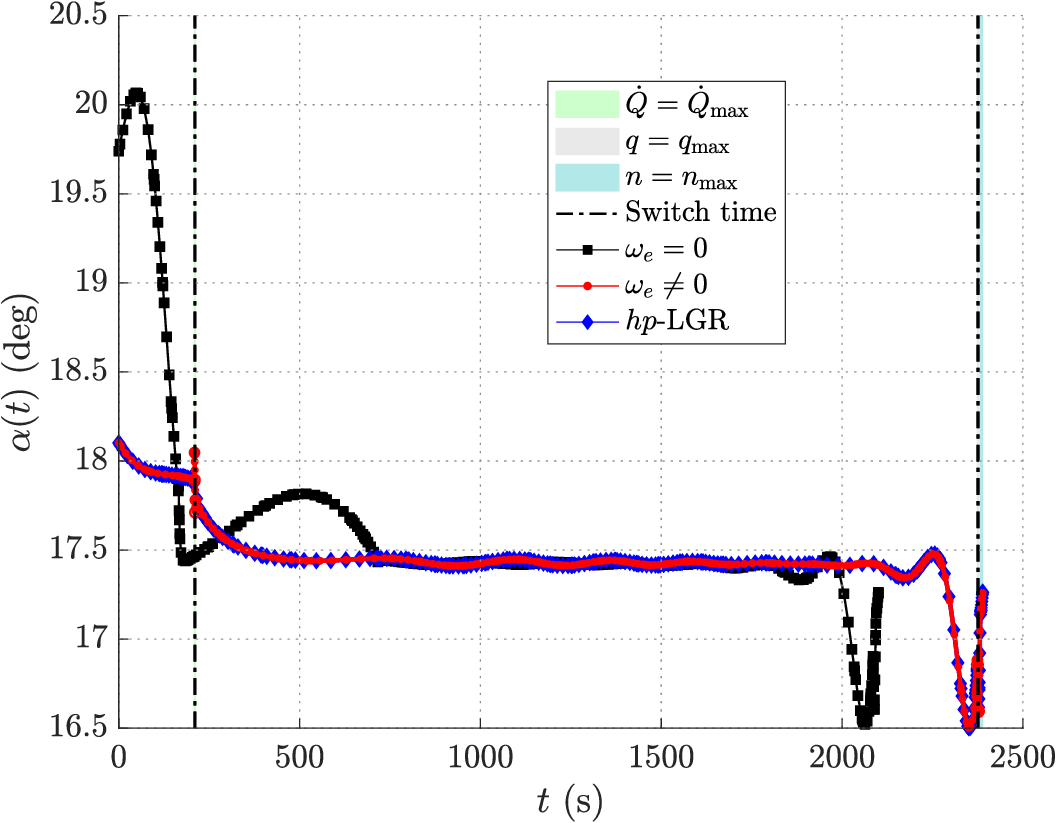}}
~~\subfloat[Bank angle, $\sigma(t)$, vs. time, $t$]{\label{fig:Study_bank}\includegraphics[scale=.4]{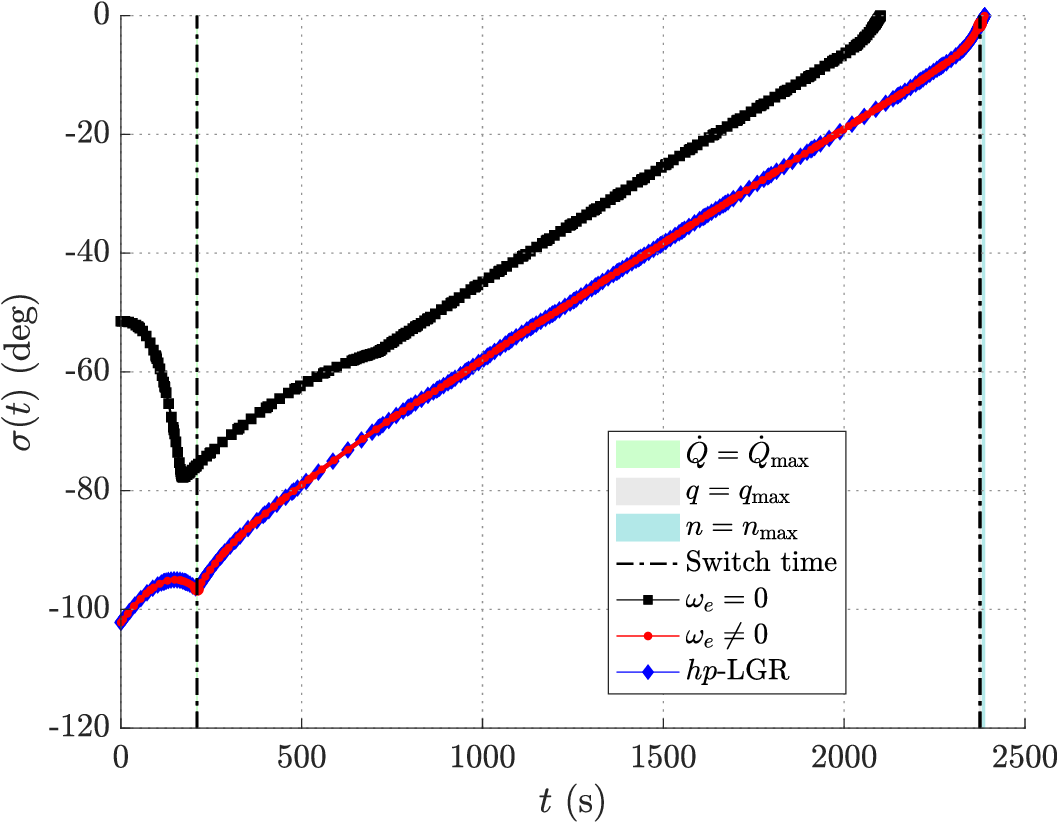}}
\caption{Optimal angle of attack, $\alpha(t)$, and bank angle, $\sigma(t)$, for Study 2 compared against the optimal angle of attack and bank angle obtained in Case 1.\label{fig:Study2_control}}
\end{figure}

Figure~\ref{fig:Study2_control} provides a comparison between the optimal angle of attack and bank angle solutions with and without Earth rotation included. It is found that the reentry vehicle starts in a steeper bank angle, gradually transitioning to a profile resembling one without considering the rotation of the Earth as it nears a level-flight configuration of $\sigma(t_f) = 0$ deg. Conversely, the vehicle begins at a shallower angle of attack than the angle of attack without Earth rotation. As a result of the change in the control structures, the control limits that were enforced in Section~\ref{section:case_2} are excluded in this study. 

From both Figs.~\ref{fig:Study2_constraints} and~\ref{fig:Study2_control} it can be seen that, in order to meet the same set of boundary conditions given in Eq.~\eqref{eq:RLVE_bc}, the total time of flight is larger when including the rotation of the Earth. Specifically, the time of flight with Earth rotation was found to be $t_{f} = 2388.62$ s while the time of flight without Earth rotation was found to be $t_{f} = 2100.47$ s. Consequently, the final geocentric latitude and Earth-relative longitude achieved is larger than the final latitude and longitude achieved without Earth rotation. Specifically, the final crossrange (latitude) achieved with Earth rotation was found to be $\phi(t_f) = 37.01$ deg, while the final crossrange achieved without Earth rotation was found to be $\phi(t_f) = 33.99$ deg. The final Earth-relative longitude achieved with Earth rotation was found to be $\theta(t_f) = 100$ deg. To compare the final Earth-relative longitude to the longitude achieved in the non-rotating Earth case, the inertial longitude is computed using the relation $\theta_I(t_f) = \theta(t_f) + \omega_e t_f$. This gives a final inertial longitude of $\theta_I(t_f) = 109.98$ deg with Earth rotation in comparison to the final longitude of $\theta_I(t_f) = 81.72$ deg achieved without Earth rotation included.

\subsubsection{Heating Rate Variation Analysis}\label{section:study_2_variation}
\begin{figure}[htb] 
\centering
\subfloat[Heating rate, $\dot{Q}(t)$, vs. time, $t$]{\label{fig:qdot_variation}\includegraphics[scale=.4]{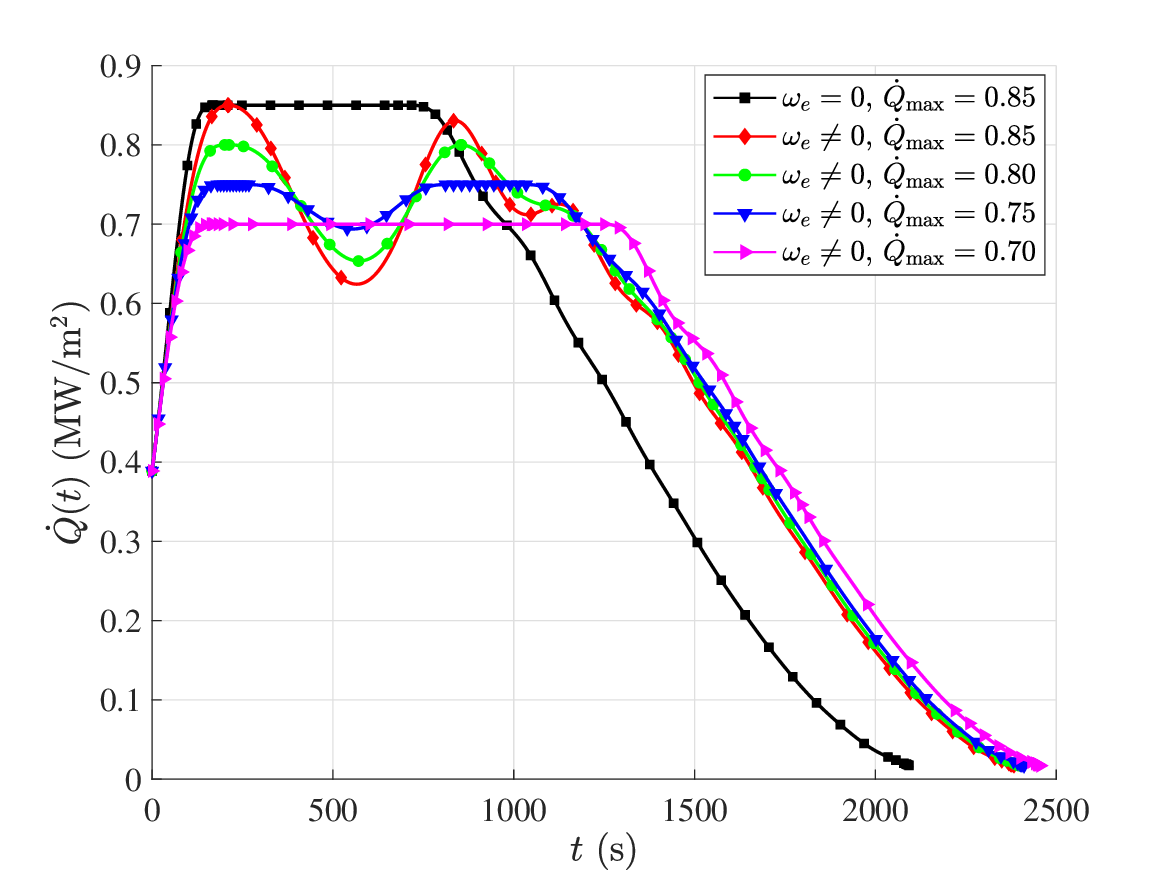}}
~~\subfloat[Dynamic pressure, $q(t)$, vs. time, $t$]{\label{fig:q_variation}\includegraphics[scale=.4]{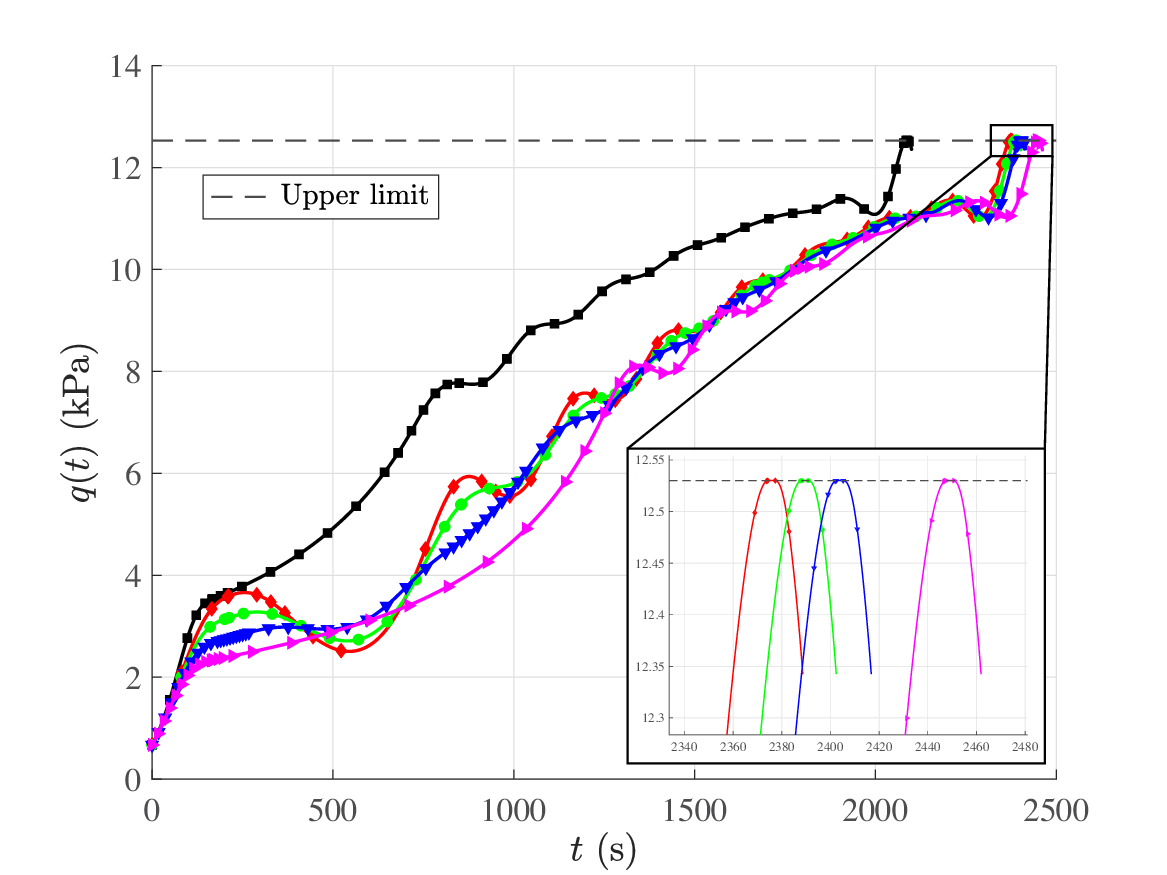}}
\caption{Optimal state-path constraint profiles for decreasing values of maximum allowable stagnation point heating rate from $\dot{Q}_{\max} = 0.85$ MW/m$^2$ to $\dot{Q}_{\max} = 0.70$ MW/m$^2$.\label{fig:Study2_statepath_variation}}
\end{figure}

In this analysis, the maximum allowable stagnation point heating rate is varied from $\dot{Q}_{\max} = 0.85$ MW/m$^2$ to $\dot{Q}_{\max} = 0.70$ MW/m$^2$. In each case, the upper limits placed on the dynamic pressure and sensed acceleration are kept at $q_{\max} = 12.53$ kPa and $n_{\max} = 1.15$, respectively. Figure~\ref{fig:Study2_statepath_variation} shows the structure of the stagnation point heating rate and dynamic pressure profiles for $\dot{Q}_{\max} = \{ 0.85, 0.80, 0.75, 0.70 \}$ MW/m$^2$ with Earth rotation included. A comparison is made against the stagnation point heating rate and dynamic pressure profiles for the solutions found in Case 1. It is seen from Fig.~\ref{fig:qdot_variation} that the $\dot{Q}_{\max}$ limit is attained twice
for $\dot{Q}_{\max} = \{0.80, 0.75 \}$ MW/m$^2$, whereas the solutions obtained for $\dot{Q}_{\max} = \{0.85, 0.70 \}$ MW/m$^2$ exhibit only one active arc along $\dot{Q}_{\max}$ of varying duration. Specifically, Table~\ref{table:optimaltimes_variation} provides the optimized entry and exit times of each active state-path constraint detected by the SPOC method. Upon further examination of Table~\ref{table:optimaltimes_variation} and Fig.~\ref{fig:q_variation}, it is seen that the occurrence of the active dynamic pressure constraint shifts with the increasing time of flight, but the duration of the constraint remains the same. Additionally, it can be found from Table~\ref{table:optimaltimes_variation} that as the limit on the stagnation point heating rate is decreased, the final latitude achieved by the reentry vehicle also decreases. Finally, it is noted that as the maximum allowable stagnation point heating rate is decreased, the time of flight increases, and consequently the total heating load on the vehicle tends to increase. Although not considered in the design of the trajectories shown, the total heating load for $\dot{Q}_{\max} = \{0.85, 0.80, 0.75, 0.70 \}$ MW/m$^2$ is found to be $Q = \{ 1244, 1249, 1252, 1250 \}$ MJ/m$^2$, respectfully. A similar trend was observed in the work of Ref.~\cite{MillerRao2021}, which studied the combined ascent-entry trajectory optimization of a common aero vehicle (CAV) model described in Ref.~\cite{Phillips2003}. 

\begin{table}[htb]
\centering
\caption{Optimized entry and exit times of each active state-path constraint for decreasing maximum allowable stagnation point heating rate.}
\begin{tabular*}{\textwidth}{@{\extracolsep{\fill}}*{8}{c}}
\hline \hline
 \rule{0pt}{3ex}   & \multicolumn{2}{c}{Arc 1, $\dot{Q} = \dot{Q}_{\max}$} & \multicolumn{2}{c}{ Arc 2, $\dot{Q} = \dot{Q}_{\max}$} & \multicolumn{2}{c}{Arc 3, $q = q_{\max}$} & Objective \\[3pt] \cline{2-3}\cline{4-5}\cline{6-7}\cline{8-8}
 \rule{0pt}{3ex} $\dot{Q}_{\max}$ (MW/m$^2$) & Entry $(s)$ & Exit $(s)$ & Entry $(s)$ & Exit $(s)$ & Entry $(s)$ & Exit $(s)$ & $\phi(t_f)$ (deg) \\[3pt] \hline
 \rule{0pt}{3ex} 0.85 & 209.45 & 211.81 & -- & -- & 2374.21 & 2377.32 & 37.01 \\
 \rule{0pt}{3ex} 0.80 & 200.63 & 212.89 & 853.83 & 855.37 & 2388.16 & 2390.89 & 36.99 \\[3pt] 
 \rule{0pt}{3ex} 0.75 & 180.77 & 267.73 & 811.23 & 1034.10 & 2402.41 & 2405.23 & 36.95 \\[3pt] 
 \rule{0pt}{3ex} 0.70 & 162.35 & 1250.24 & -- & -- & 2447.47 & 2450.85 & 36.84 \\[3pt] 
\hline \hline
\end{tabular*}
\label{table:optimaltimes_variation}
\end{table}

\begin{figure}[htb] 
\vspace{-0.5em}
\centering
\subfloat[Angle of attack, $\alpha(t)$, vs. time, $t$]{\label{fig:Study2_aoa_variation}\includegraphics[scale=.4]{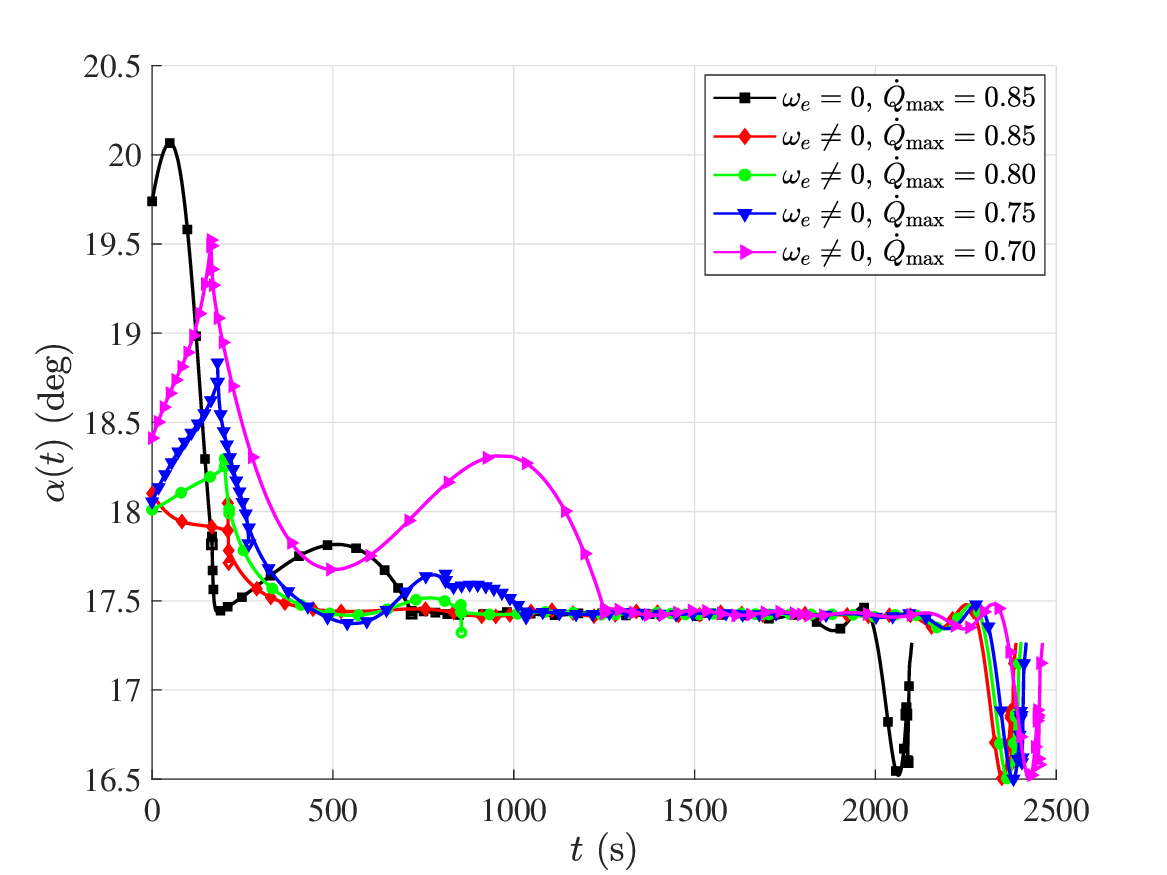}}
~~\subfloat[Bank angle, $\sigma(t)$, vs. time, $t$]{\label{fig:Study2_bank_variation}\includegraphics[scale=.4]{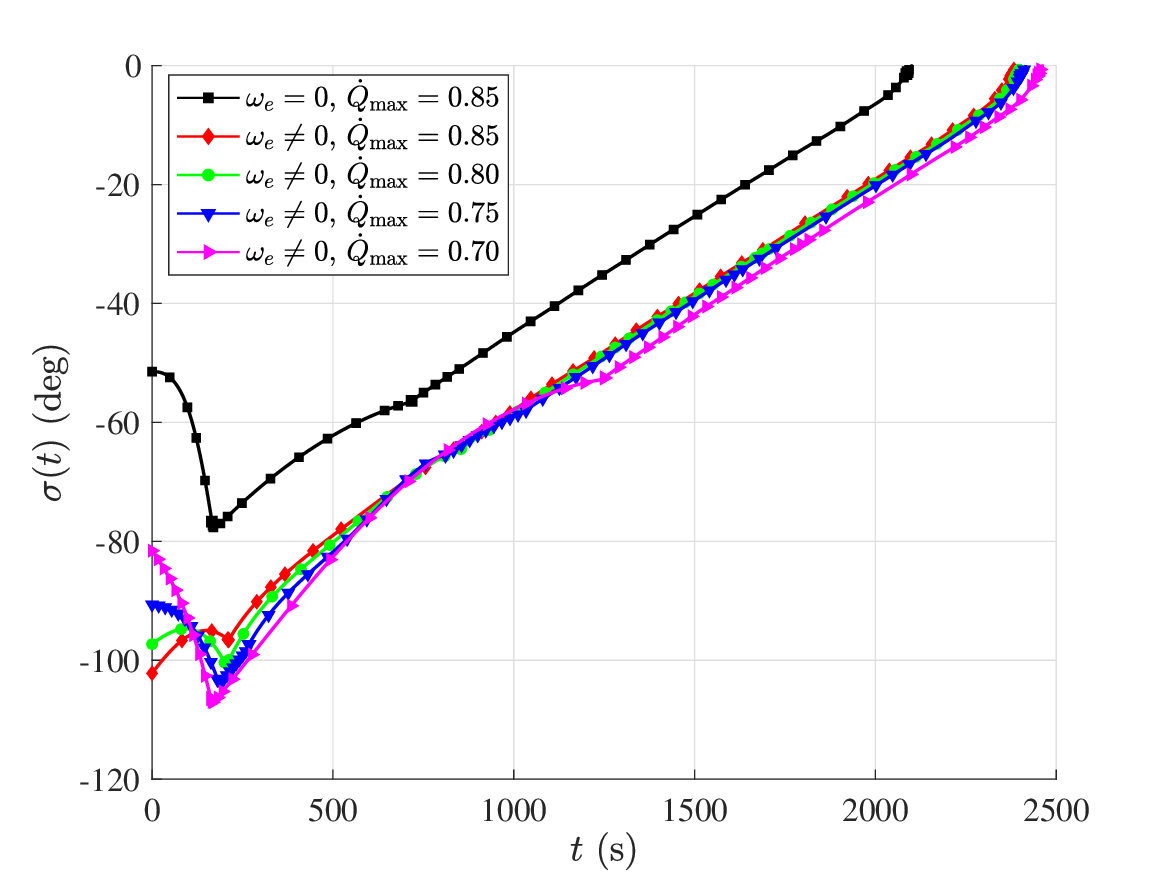}}
\caption{Optimal angle of attack, $\alpha(t)$, and bank angle, $\sigma(t)$, for decreasing values of maximum allowable stagnation point heating rate from $\dot{Q}_{\max} = 0.85$ MW/m$^2$ to $\dot{Q}_{\max} = 0.70$ MW/m$^2$.\label{fig:Study2_control_variation}}
\end{figure}

Figure~\ref{fig:Study2_control_variation} shows the optimal angle of attack and bank angle for $\dot{Q}_{\max} = \{ 0.85, 0.80, 0.75, 0.70 \}$ MW/m$^2$ with Earth rotation included. A comparison is also made against the angle of attack and bank angle profiles found in Case 1. It can be seen in Fig.~\ref{fig:Study2_aoa_variation} that across all cases the angle of attack spikes when encountering the upper limit of the stagnation point heating rate, where the magnitude of the spike increases with decreasing values of $\dot{Q}_{\max}$. Next, it can observed in Fig.~\ref{fig:Study2_bank_variation} that for decreasing values of $\dot{Q}_{\max}$, the vehicle starts in a shallower bank angle and then achieves a steeper bank angle as it approaches $\dot{Q}_{\max}$. After achieving the steepest bank angle, each case follows a similar profile for the remainder of the trajectory. In every case, the bank angle is steeper with Earth rotation included versus without Earth rotation, whereas the angle of attack generally attains a similar profile for both problem formulations.

\vspace{-2em}
\subsection{Discussion on Control Discontinuities}\label{section:control_discontinuity}
\begin{figure}[htb!] 
\centering
\includegraphics[scale=.4]{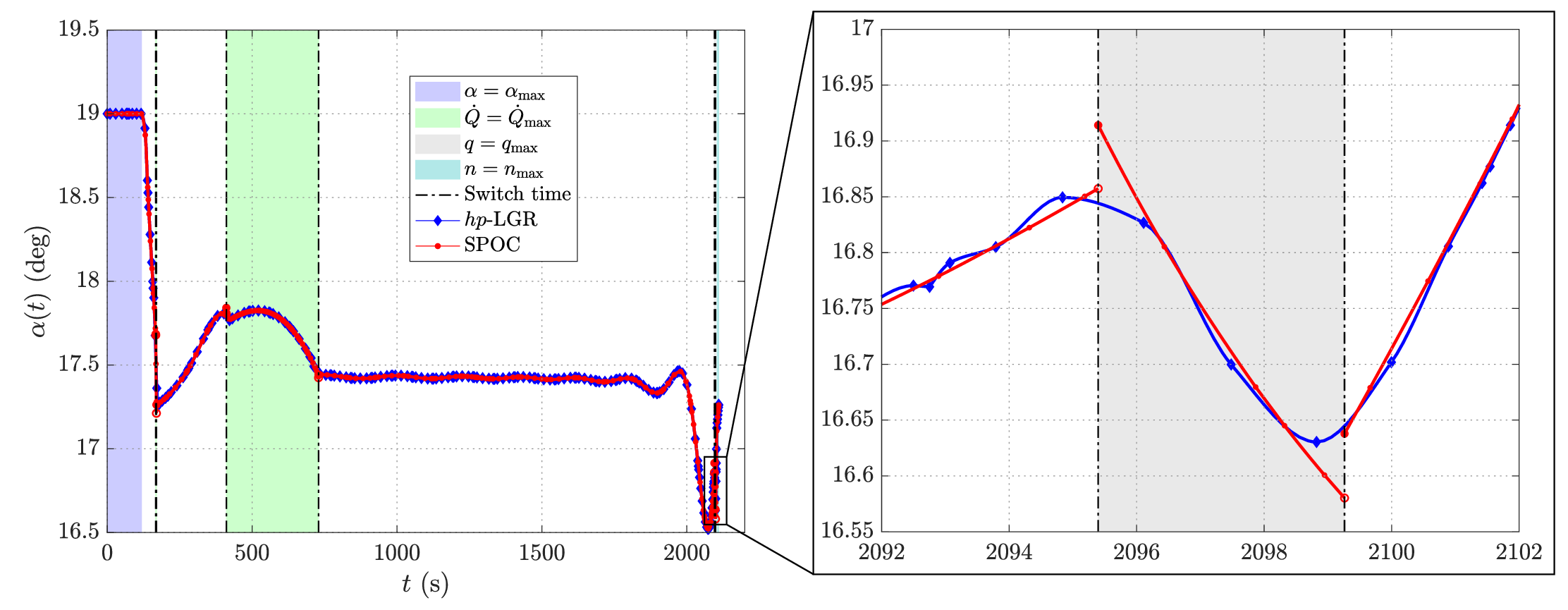}
\caption{Small discontinuities captured by the SPOC method for the optimal angle of attack along the active dynamic pressure constraint arc for Case 2.\label{fig:Case2_aoa_qarc}}
\end{figure}

The SPOC method captures small discontinuities in the angle of attack at the entry and exit times of each active state-path constraint. For example, Fig.~\ref{fig:Case2_aoa_qarc} provides an enlarged view of the angle of attack computed by the SPOC method for Case 2 along the active dynamic pressure constraint. Because the multiple-domain LGR formulation does not include a continuity constraint on the angle of attack at an interface between two domains, it is possible for the angle of attack to be discontinuous at a domain interface time.    
Moreover, examining the form of Eqs.~\eqref{eq:vdot}, ~\eqref{eq:heatrate}, and~\eqref{eq:dynamicpressure} for the rate of change of the speed, the heating rate, and the dynamic pressure, respectively, the angle of attack is the control that appears explicitly in the first time derivative of both state-path constraints. Consequently, Eq.~\eqref{eq:mixedconstr} is a constraint that depends upon the angle of attack and is enforced on the active state-path constraint arcs. 
This study shows that by implementing the multiple-domain LGR formulation, the SPOC method (i) captures control discontinuities occurring at the entry and exit of active state-path constraints, and (ii) computes separate control arcs within each domain to obtain an accurate approximation of the active state-path constraints.

\section{Conclusions\label{section:conclusions}}
The trajectory optimization problem of maximizing the crossrange of a reusable launch vehicle subject to two state-inequality path constraints, two control-inequality path constraints, and one mixed state-control inequality path constraint is studied. A recently developed method for solving state-path constrained optimal control (SPOC) problems is used. The SPOC method 
automatically detects the state-path constraint structure and partitions the domain of the independent variable into multiple subdomains where the state-path constraints are either active or inactive. Additional necessary conditions are then algorithmically enforced within the subdomains that contain active state-path constraints, and the detected entry and exit times are included as additional decision variables in the optimization. Two studies were performed which analyzed a variety of problem formulations. In the first study, the SPOC method was used to generate results for two cases where the control-inequality path constraints are 1) inactive and 2) active, respectively. It was found that by enforcing the control constraints, the optimal heating rate constraint profile changes from containing a single active constraint arc to having two active constraint arcs, respectively. Next, the SPOC method was used to generate new results to the reusable launch vehicle optimal control problem with Earth rotational effects. It was found that with the inclusion of the rotation of the Earth, the total time of flight increases and the peak heating rate on the vehicle decreases, leading to a change in the structure of the active heating rate constraint. Finally, a feature of the SPOC method to compute separate control arcs within each domain and capture small discontinuities in the control is discussed. In all of the problem formulations studied, it is found that, with minimal user input, the SPOC method is able to detect and optimize changes in the active constraint structure as well as compute accurate solutions to a complex constrained hypersonic reentry problem. 

\section{Acknowledgments}
The authors gratefully acknowledge support for this research from the U.S.~National Science Foundation under grant CMMI-2031213, the U.S.~Office of Naval Research under grant N00014-22-1-2397, and from the U.S.~Air Force Research Laboratory under contract FA8651-21-F-1041.  




\renewcommand{\baselinestretch}{1.0}
\normalsize\normalfont

\begin{thebibliography}{10}
\newcommand{\enquote}[1]{``#1''}

\bibitem{Dickmanns1972}
Dickmanns, E.~D., \enquote{Maximum {R}ange {T}hree-{D}imensional {L}ifting
  {P}lanetary {E}ntry,} Tech. rep., National Aeronautics and Space
  Administration, 1972.

\bibitem{ZondervanHuffman1984}
Zondervan, K., Bauer, T., Betts, J., and Huffman, W., \enquote{{S}olving the
  {O}ptimal {C}ontrol {P}roblem {U}sing a {N}onlinear {P}rogramming {T}echnique
  {P}art 3: {O}ptimal {S}huttle {R}eentry {T}rajectories,} {\em Proceedings of
  the Astrodynamics Specialist Conference\/}, AIAA-84-2039, Seattle, WA, August
  20-22, 1984, pp. 1--8.

\bibitem{KugelmannPesch1990}
Kugelmann, B. and Pesch, H.~J., \enquote{New general guidance method in
  constrained optimal control, Part 2: Application to space shuttle guidance,}
  {\em Journal of Optimization Theory and Applications\/}, Vol.~67, No.~3,
  1990, pp.~437--446.

\bibitem{Pesch1994}
Pesch, H.~J., \enquote{A practical guide to the solution of real-life optimal
  control problems,} {\em Control and cybernetics\/}, Vol.~23, No.~1, 1994,
  pp.~2.

\bibitem{Betts2020}
Betts, J.~T., {\em Practical Methods for Optimal Control Using Nonlinear
  Programming\/}, SIAM Press, Philadelphia, PA, 2020.

\bibitem{DarbyRao2011b}
Darby, C.~L., Hager, W.~W., and Rao, A.~V., \enquote{{D}irect {T}rajectory
  {O}ptimization {U}sing a {V}ariable {L}ow-{O}rder {A}daptive {P}seudospectral
  {M}ethod,} {\em Journal of Spacecraft and Rockets\/}, Vol.~48, No.~3,
  May--June 2011, pp.~433--445.

\bibitem{PattersonRao2015}
Patterson, M.~A., Hager, W.~W., and Rao, A.~V., \enquote{{A} $ph$ {M}esh
  {R}efinement {M}ethod for {O}ptimal {C}ontrol,} {\em Optimal Control
  Applications and Methods\/}, Vol.~36, No.~4, July--August 2015, pp.~398--421.

\bibitem{DennisRao2019}
Dennis, M.~E., Hager, W.~W., and Rao, A.~V., \enquote{{C}omputational {M}ethod
  for {O}ptimal {G}uidance and {C}ontrol {U}sing {A}daptive {G}aussian
  {Q}uadrature {C}ollocation,} {\em Journal of Guidance, Control, and
  Dynamics\/}, Vol.~42, No.~9, September 2019, pp.~2026--2041.

\bibitem{MallTaheri2022}
Mall, K. and Taheri, E., \enquote{Three-Degree-of-Freedom Hypersonic Reentry
  Trajectory Optimization Using an Advanced Indirect Method,} {\em Journal of
  Spacecraft and Rockets\/}, Vol.~59, No.~5, 2022, pp.~1463--1474.

\bibitem{RaoTang2002}
Rao, A.~V., Tang, S., and Hallman, W.~P., \enquote{Numerical optimization study
  of multiple-pass aeroassisted orbital transfer,} {\em Optimal Control
  Applications and Methods\/}, Vol.~23, No.~4, May-June 2002, pp.~215--238.

\bibitem{DarbyRao2011c}
Darby, C.~L. and Rao, A.~V., \enquote{Minimum-{F}uel {L}ow-{E}arth-{O}rbit
  {A}eroassisted {O}rbital {T}ransfer of {S}mall {S}pacecraft,} {\em Journal of
  Spacecraft and Rockets\/}, July-August 2011.

\bibitem{FuhrRao2018}
Fuhr, R. and Rao, A.~V., \enquote{Minimum-{F}uel {L}ow-{E}arth-{O}rbit
  {A}eroglide and {A}erothrust {A}eroassisted {O}rbital {T}ransfer {S}ubject to
  {H}eating {C}onstraints,} {\em Journal of Spacecraft and Rockets\/}, May-June
  2018.

\bibitem{GrantBraun2015}
Grant, M.~J. and Braun, R.~D., \enquote{Rapid indirect trajectory optimization
  for conceptual design of hypersonic missions,} {\em Journal of Spacecraft and
  Rockets\/}, Vol.~52, No.~1, 2015, pp.~177--182.

\bibitem{ZhengAi2017}
Zheng, Y., Cui, H., and Ai, Y., \enquote{Indirect trajectory optimization for
  mars entry with maximum terminal altitude,} {\em Journal of Spacecraft and
  Rockets\/}, Vol.~54, No.~5, 2017, pp.~1068--1080.

\bibitem{MillerRao2021}
Miller, A.~T. and Rao, A.~V., \enquote{End-to-{E}nd {P}erformance
  {O}ptimization for {H}igh-{S}peed {A}scent {E}ntry {M}issions,} {\em Journal
  of Spacecraft and Rockets\/}, December 2021.

\bibitem{BrysonHo1975}
Bryson, A.~E. and Ho, Y., {\em Applied Optimal Control: Optimization,
  Estimation, and Control\/}, Hemisphere Publishing Corporation, 1975.

\bibitem{Betts1998}
Betts, J.~T., \enquote{{S}urvey of {N}umerical {M}ethods for {T}rajectory
  {O}ptimization,} {\em Journal of Guidance, Control, and Dynamics\/}, Vol.~21,
  No.~2, March-April 1998, pp.~193--207.

\bibitem{BertrandEpenoy2002}
Bertrand, R. and Epenoy, R., \enquote{New smoothing techniques for solving
  bang--bang optimal control problems—numerical results and statistical
  interpretation,} {\em Optimal Control Applications and Methods\/}, Vol.~23,
  No.~4, 2002, pp.~171--197.

\bibitem{GraichenChaplais2010}
Graichen, K., Kugi, A., Petit, N., and Chaplais, F., \enquote{Handling
  constraints in optimal control with saturation functions and system
  extension,} {\em Systems \& Control Letters\/}, Vol.~59, No.~11, 2010,
  pp.~671--679.

\bibitem{Epenoy2011}
Epenoy, R., \enquote{Fuel optimization for continuous-thrust orbital rendezvous
  with collision avoidance constraint,} {\em Journal of Guidance, Control, and
  Dynamics\/}, Vol.~34, No.~2, 2011, pp.~493--503.

\bibitem{TaheriAtkins2016}
Taheri, E., Kolmanovsky, I., and Atkins, E., \enquote{Enhanced smoothing
  technique for indirect optimization of minimum-fuel low-thrust trajectories,}
  {\em Journal of Guidance, Control, and Dynamics\/}, Vol.~39, No.~11, 2016,
  pp.~2500--2511.

\bibitem{AntonyGrant2018}
Antony, T. and Grant, M.~J., \enquote{Path constraint regularization in optimal
  control problems using saturation functions,} {\em 2018 AIAA Atmospheric
  Flight Mechanics Conference\/}, 2018, p. 0018.

\bibitem{MansellGrant2018}
Mansell, J.~R. and Grant, M.~J., \enquote{Adaptive continuation strategy for
  indirect hypersonic trajectory optimization,} {\em Journal of Spacecraft and
  Rockets\/}, Vol.~55, No.~4, 2018, pp.~818--828.

\bibitem{MallTaheri2020}
Mall, K., Grant, M.~J., and Taheri, E., \enquote{Uniform Trigonometrization
  Method for Optimal Control Problems with Control and State Constraints,} {\em
  Journal of Spacecraft and Rockets\/}, Vol.~57, No.~5, sep 2020,
  pp.~995--1007.

\bibitem{HeidrichGrant2021}
Heidrich, C.~R., Sparapany, M.~J., and Grant, M.~J., \enquote{Generalized
  Regularization of Constrained Optimal Control Problems,} {\em Journal of
  Spacecraft and Rockets\/}, jan 2021.

\bibitem{BensonRao2006}
Benson, D.~A., Huntington, G.~T., Thorvaldsen, T.~P., and Rao, A.~V.,
  \enquote{{D}irect {T}rajectory {O}ptimization and {C}ostate {E}stimation via
  an {O}rthogonal {C}ollocation {M}ethod,} {\em Journal of Guidance, Control,
  and Dynamics\/}, Vol.~29, No.~6, November-December 2006, pp.~1435--1440.

\bibitem{KameswaranBiegler2008a}
Kameswaran, S. and Biegler, L.~T., \enquote{{C}onvergence {R}ates for {D}irect
  {T}ranscription of {O}ptimal {C}ontrol {P}roblems {U}sing {C}ollocation at
  {R}adau {P}oints,} {\em Computational Optimization and Applications\/},
  Vol.~41, No.~1, 2008, pp.~81--126.

\bibitem{GargHuntington2010}
Garg, D., Patterson, M.~A., Hager, W.~W., Rao, A.~V., Benson, D.~A., and
  Huntington, G.~T., \enquote{{A} {U}nified {F}ramework for the {N}umerical
  {S}olution of {O}ptimal {C}ontrol {P}roblems {U}sing {P}seudospectral
  {M}ethods,} {\em Automatica\/}, Vol.~46, No.~11, November 2010,
  pp.~1843--1851. DOI: 10.1016/j.automatica.2010.06.048.

\bibitem{GargRao2011a}
Garg, D., Hager, W.~W., and Rao, A.~V., \enquote{{P}seudospectral {M}ethods for
  {S}olving {I}nfinite-{H}orizon {O}ptimal {C}ontrol {P}roblems,} {\em
  Automatica\/}, Vol.~47, No.~4, April 2011, pp.~829--837. DOI:
  10.1016/j.automatica.2011.01.085.

\bibitem{GargRao2011b}
Garg, D., Patterson, M.~A., Darby, C.~L., Francolin, C., Huntington, G.~T.,
  Hager, W.~W., and Rao, A.~V., \enquote{{D}irect {T}rajectory {O}ptimization
  and {C}ostate {E}stimation of {F}inite-{H}orizon and {I}nfinite-{H}orizon
  {O}ptimal {C}ontrol {P}roblems via a {R}adau {P}seudospectral {M}ethod,} {\em
  Computational Optimization and Applications\/}, Vol.~49, No.~2, June 2011,
  pp.~335--358. DOI: 10.1007/s10589--00--09291--0.

\bibitem{HagerWang2018}
Hager, W.~W., Liu, J., Mohapatra, S., Rao, A.~V., and Wang, X.-S.,
  \enquote{Convergence Rate for a Gauss Collocation Method Applied to
  Constrained Optimal Control,} {\em SIAM Journal on Control and
  Optimization\/}, Vol.~56, No.~2, 2018, pp.~1386--1411.

\bibitem{Elnagar1995}
Elnagar, G., Kazemi, M.~A., and Razzaghi, M., \enquote{The {P}seudospectral
  {L}egendre {M}ethod for {D}iscretizing {O}ptimal {C}ontrol {P}roblems,} {\em
  IEEE transactions on Automatic Control\/}, Vol.~40, No.~10, 1995,
  pp.~1793--1796.

\bibitem{HouRao2012}
Hou, H., Hager, W., and Rao, A., \enquote{Convergence of a Gauss Pseudospectral
  Method for Optimal Control,} {\em {AIAA} Guidance, Navigation, and Control
  Conference\/}, American Institute of Aeronautics and Astronautics, aug 2012,
  p. 4452.

\bibitem{DarbyRao2011a}
Darby, C.~L., Hager, W.~W., and Rao, A.~V., \enquote{An $hp$‐-{A}daptive
  {P}seudospectral {M}ethod for {S}olving {O}ptimal {C}ontrol {P}roblems,} {\em
  Optimal Control Applications and Methods\/}, Vol.~32, No.~4, 2011,
  pp.~476--502.

\bibitem{LiuRao2015}
Liu, F., Hager, W.~W., and Rao, A.~V., \enquote{Adaptive mesh refinement method
  for optimal control using nonsmoothness detection and mesh size reduction,}
  {\em Journal of the Franklin Institute\/}, Vol.~352, No.~10, oct 2015,
  pp.~4081--4106.

\bibitem{LiuRao2018}
Liu, F., Hager, W.~W., and Rao, A.~V., \enquote{{A}daptive {M}esh {R}efinement
  {M}ethod for {O}ptimal {C}ontrol {U}sing {D}ecay {R}ates of {L}egendre
  {P}olynomial {C}oefficients,} {\em IEEE Transactions on Control System
  Technology\/}, Vol.~26, No.~4, July 2018, pp.~1475--1483.

\bibitem{PagerRao2022}
Pager, E.~R. and Rao, A.~V., \enquote{Method for solving bang-bang and singular
  optimal control problems using adaptive Radau collocation,} {\em
  Computational Optimization and Applications\/}, Vol.~59, No.~81, Jan 2022,
  pp.~857--887.

\bibitem{ByczkowskiRao2024a}
Byczkowski, C.~A. and Rao, A.~V., \enquote{Method for solving state-path
  constrained optimal control problems using adaptive Radau collocation,} {\em
  Optimal Control Applications and Methods\/}, Vol.~45, No.~3, May 2024,
  pp.~1199--1222.

\bibitem{FeeheryBarton1998}
Feehery, W.~F. and Barton, P.~I., \enquote{Dynamic optimization with state
  variable path constraints,} {\em Computers \& Chemical Engineering\/},
  Vol.~22, No.~9, 1998, pp.~1241--1256.

\bibitem{ByczkowskiRao2023a}
Byczkowski, C.~A. and Rao, A.~V., \enquote{{S}tructure {D}etection {M}ethod for
  {S}olving {S}tate {V}ariable {I}nequality {P}ath {C}onstrained {O}ptimal
  {C}ontrol {P}roblems,} {\em Proceedings of the AAS Space Flight Mechanics
  Meeting\/}, AAS Paper 23-165, Austin, TX, January 15-19, 2023, p. 0013.

\bibitem{ByczkowskiRao2024b}
Byczkowski, C.~A. and Rao, A., \enquote{Constrained Hypersonic Reentry
  Trajectory Optimization Using A Multiple-Domain Direct Collocation Method,}
  {\em AIAA SCITECH 2024 Forum\/}, 2024, p. 1457.

\bibitem{BrysonDreyfus1963}
Bryson, A.~E., Denham, W.~F., and Dreyfus, S.~E., \enquote{Optimal Programming
  Problems with Inequality Constraints I: Necessary Conditions for Extremal
  Solutions,} {\em AIAA Journal\/}, Vol.~1, No.~11, November 1963,
  pp.~2544--2550.

\bibitem{SuttonGraves1971}
Sutton, K. and Graves~Jr, R.~A., \enquote{A General Stagnation-Point Convective
  Heating Equation For Arbitrary Gas Mixtures,} Tech. Rep. TR R-376, NASA
  Langley Research Center, 1971.

\bibitem{PattersonRao2014}
Patterson, M.~A. and Rao, A.~V., \enquote{{GPOPS}-{II}: {A} {MATLAB} {S}oftware
  for {S}olving {M}ultiple-{P}hase {O}ptimal {C}ontrol {P}roblems {U}sing
  $hp$-{A}daptive {G}aussian {Q}uadrature {C}ollocation {M}ethods and {S}parse
  {N}onlinear {P}rogramming,} {\em {ACM} Transactions on Mathematical
  Software\/}, Vol.~41, No.~1, oct 2014, pp.~1--37.

\bibitem{BieglerZavala2008}
Biegler, L.~T. and Zavala, V.~M., \enquote{{L}arge-{S}cale {N}onlinear
  {P}rogramming {U}sing {IPOPT}: {A}n {I}ntegrating {F}ramework for
  {E}nterprise-{W}ide {O}ptimization,} {\em Computers and Chemical
  Engineering\/}, Vol.~33, No.~3, March 2008, pp.~575--582.

\bibitem{ma57}
Duff, I.~S., \enquote{MA57---a Code for the Solution of Sparse Symmetric
  Definite and Indefinite Systems,} {\em ACM Trans. Math. Softw.\/}, Vol.~30,
  No.~2, June 2004, pp.~118–144.

\bibitem{PattersonRao2012}
Patterson, M.~A. and Rao, A.~V., \enquote{Exploiting sparsity in direct
  collocation pseudospectral methods for solving optimal control problems,}
  {\em Journal of Spacecraft and Rockets\/}, Vol.~49, No.~2, 2012,
  pp.~354--377.

\bibitem{ByczkowskiRao2023b}
Byczkowski, C.~A. and Rao, A.~V., \enquote{{R}eusable {E}ntry {V}ehicle
  {T}rajectory {O}ptimization {U}sing {M}ultiple-{D}omain {R}adau
  {C}ollocation,} {\em Proceedings of the AIAA SCITECH 2023 Forum\/}, AIAA
  Paper 2023-1168, National Harbor, MD, January 23-27, 2023, p. 1168.

\bibitem{VinhCulp1980}
Vinh, N.-X., Busemann, A., and Culp, R.~D., {\em Hypersonic and Planetary Entry
  Flight Mechanics\/}, University of Michigan Press, Ann Arbor, Michigan, 1980.

\bibitem{Phillips2003}
Phillips, T.~H., \enquote{A common aero vehicle (CAV) model, description, and
  employment guide,} {\em Schafer Corporation for AFRL and AFSPC\/}, Vol.~27,
  2003, pp.~1--12.

\end{thebibliography}

\end{document}